\documentclass{art_ar}
\usepackage{amssymb}
\usepackage{amsfonts}
\usepackage{array}
\usepackage{longtable}
\usepackage{tabularx}

\usepackage{makeidx}
\usepackage{amsmath}
\usepackage{theorem}
\usepackage[mathscr]{eucal}
\usepackage{enumerate}
\usepackage[dvips]{epsfig}
\usepackage{float}  
\usepackage{subfigure} 
\usepackage{color}       

\theoremstyle{plain}
\newtheorem{Def}{Definition}[section]

\newtheorem{The}[Def]{Theorem}

\newtheorem{Bem}[Def]{Remark}
\newtheorem{Lem}[Def]{Lemma}

\numberwithin{equation}{section}

\newcommand{\E}{\operatorname{E}}
\newcommand{\Prob}{\operatorname{P}}

\begin{document}

\begin{frontmatter}

\title{Stability Analysis and Classification of Runge-Kutta Methods for
Index~1 Stochastic Differential-Algebraic Equations \\
with Scalar Noise}
\author{Dominique K{\"u}pper}, 
\author{Anne Kv\ae{}rn\o{}\thanksref{label2}},
\author{Andreas R{\"o}{\ss}ler\thanksref{label3}}
\ead{roessler@math.uni-luebeck.de}
%
\address[label2]{Department of Mathematical Sciences, Norwegian University of
Science and Technology, Norway}
\address[label3]{Universit\"at zu L\"ubeck, Institut f\"ur Mathematik,
Ratzeburger Allee 160, D-23562~L\"ubeck, Germany}




%
\begin{abstract}
    The problem of solving stochastic differential-algebraic
    equations (SDAEs)
    of index one with a scalar driving Brownian motion is considered.
    Recently, the authors proposed a class of stiffly accurate
    stochastic Runge-Kutta (SRK) methods that do not involve any pseudo-inverses
    or projectors for the numerical solution of the problem. Based on this class of
    approximation methods, a classification for the coefficients of
    stiffly accurate SRK methods attaining strong order 0.5 as well as strong order 1.0
    are calculated. Further, the mean-square stability for
    the considered class of SRK methods is analysed. As the main result, families of
    $A$-stable efficient order 0.5 and 1.0 stiffly accurate SRK methods with a
    minimal number of stages for SDEs as well as for SDAEs are presented.
\end{abstract}
\end{frontmatter}
\section{Introduction}
In many applications like, e.~g., the simulation of the dynamics of
multibody systems, optimal control problems or electric circuit
simulation (see \cite{Asch,LoTh12,SWW09} for more details),
differential-algebraic equations serve as a model for the dynamical
system under consideration. However, often random disturbances, that
can be described by some noise source, have to be taken into
account. This leads to models based on stochastic
differential-algebraic equations (SDAEs) and numerical solutions
need to be calculated whenever explicit solutions are not available.
In~\cite{KKR12}, the authors propose a class of stiffly accurate
stochastic Runge-Kutta (SRK) methods that can be applied for the
numerical solution of nonlinear index~1 SDAEs with scalar noise. The
introduced class of SRK methods contains schemes attaining orders of
convergence 0.5 and 1.0 in the mean-square sense. Compared to well
known numerical schemes for SDAEs (see \cite{KKR12} for details),
their main advantages are that they do not need the calculation of
any pseudo-inverses or
projectors and can be applied directly to the SDAE system. \\ \\
In the following, we first give a classification of the space of
solutions for order 0.5 and order 1.0 conditions derived in
\cite{KKR12} in case of stiffly accurate methods that are diagonally
implicit in the drift part. Based on this classification, we
determine some coefficients for the SRK method such that the number
of stages is minimal in order to reduce computational costs.
Applying the calculated classification yields the main result: We
present families of stiffly accurate SRK methods for which
$A$-stability is proven explicitly and that have a minimal
number of stages and implicit equations to be solved each step. \\ \\
The paper is organized as follows: In
Section~\ref{Section:Class-SRK-SDAEs}, we present the general class
of SRK methods under consideration, that can be applied to index~1
SDAE systems with scalar noise. Especially, the strong order
conditions for the SRK methods calculated in \cite{KKR12} are given,
representing the basis for the classification of order 0.5 SRK
methods in Section~\ref{Sec:Class-order-0.5} and of order 1.0 SRK
methods in Section~\ref{Sec:Class-order-1.0}. The classification is
then used in Section~\ref{Sec:Optimal-Scheme-Stability-Ana} in order
to determine some coefficients for schemes with a minimal number of
stages and to analyse their mean-square stability properties.
Finally, 
some families of $A$-stable SRK methods are presented and their
$A$-stability is proved explicitly.
\section{Stochastic Runge-Kutta Methods for SDAEs}
\label{Section:Class-SRK-SDAEs}
Let $(\Omega, \mathcal{F}, \Prob)$ be a complete probability space
equipped with a filtration $(\mathcal{F}_t)_{t \geq 0}$ fulfilling
the usual conditions. Further, let $(W_t)_{t \geq 0}$ be a real
valued Wiener process adapted to $(\mathcal{F}_t)_{t \geq 0}$ and
let $\mathcal{I}=[t_0,T]$ for some $0 \leq t_0 < T$. Then, we denote
by $(X_t)_{t \in \mathcal{I}}$ the $d$-dimensional solution of the
index~1 It{\^o} stochastic differential-algebraic equation system
\begin{equation}\label{SDAE-short}
    M \, {\mathrm{d}}X_t = f(t,X_t) \, {\mathrm{d}}t + g(t,X_t)
    \, {\mathrm{d}} W_t
\end{equation}
with consistent initial value $X_{t_0} \in L^2(\Omega)$. Here, $f, g
: \mathcal{I} \times \mathbb{R}^d \to \mathbb{R}^d$ are assumed to
be globally Lipschitz continuous functions and $M \in \mathbb{R}^{d
\times d}$ is a matrix. If $M$ is non-singular, multiplying by
$M^{-1}$ transforms \eqref{SDAE-short} to a classical system of
stochastic differential equations (SDEs). However, if $M$ is
singular, we have a system of SDAEs that can be written as a system
of SDEs with some algebraic constraints, see e.~g.\ \cite{KKR12}. In
this case, we assume that the noise sources do not appear in the
algebraic constraints and that the constraints are globally uniquely
solvable for algebraic variables. This guarantees that
\eqref{SDAE-short} is an index~1 SDAE system \cite{KKR12,Wi}. In the
following, we always assume that the a unique solution of
\eqref{SDAE-short} exists, see \cite{Wi} for details. Because $f$
and $g$ need not to be linear, we are concerned with a general
nonlinear system of
index~1 SDAEs driven by a scalar Wiener process. \\ \\
In order to solve \eqref{SDAE-short} numerically, we consider the
class of stiffly accurate SRK methods for the strong approximation
of $(X_t)_{t \in \mathcal{I}}$ introduced in \cite{KKR12}. The
advantage of stiffly accurate SRK methods is that they can be
directly applied to the index~1 SDAE system~\eqref{SDAE-short}.
We consider a discretization
$\mathcal{I}_h=\{t_0, t_1, \ldots t_N\}$ of $\mathcal{I}$ and we
denote by $y_n$ the approximation of $(X_t)_{t \in \mathcal{I}}$ at
time $t_n$ using step sizes $h_n = t_{n+1}-t_n>0$. Further, let
$I_{(1),n} = W_{t_{n+1}}-W_{t_n}$ denote an increment of the Wiener
process and let $I_{(1,1),n} = \frac{1}{2}(I_{(1),n}^2-h_n)$ denote
the corresponding double integral.
Then, the approximations calculated by a stiffly accurate $s$-stages
SRK method are defined by $y_0=X_{t_0}$ and
\begin{equation} \label{SSRK}
  \begin{split}
    M \cdot H_i &= {} M \cdot y_n + \sum_{j=1}^s A_{ij} h_n \, f(t_n+c_j
    h_n,H_j) \\
    & {} \quad + \sum_{j=1}^s \Big( B_{ij}^{(1)} I_{(1),n} + B_{ij}^{(2)}
    \frac{I_{(1,1),n}}{\sqrt{h_n}} + B_{ij}^{(3)} \sqrt{h_n} \Big) \, g(t_n+c_j
    h_n,H_j), \\
    y_{n+1} &= {} H_s
  \end{split}
\end{equation}
for $i = 1, \ldots, s$ and $n=0,1,\ldots,N-1$, provided that the
coefficient matrix $A=(A_{ij})$ is nonsingular or provided that the
first stage of the method is explicit with $M \cdot H_1=M \cdot y_n$
and $(A_{ij})_{i,j=2}^s$ is nonsingular, see also \cite{Kue,KKR12}.
In general, a SRK method for SDEs, see e. g.\ \cite{Roes}, is called
stiffly accurate if its last stage coincides with the approximation
rule, i. e., if $y_{n+1}=H_s$.
The SRK method \eqref{SSRK} with $s$ stages is defined by its
coefficients $A=(A_{ij})$, $B^{(k)}=(B^{(k)}_{ij})$ for $k=1,2,3$
and $c=(c_j)$ for $i,j=1, \ldots, s$ that are usually given by an
extended Butcher tableau:
\renewcommand{\arraystretch}{1.6}
\begin{equation} \label{Butcher-tableau-Ito-ord1.0-scalarnoise}
{
\begin{tabular}{r|c|c|c|c}
    $c$ & $\,\,\, A \,\,\,$ & $B^{(1)}$ & $B^{(2)}$ & $B^{(3)}$ \\
    \hline
\end{tabular}
}
\end{equation}
In order to analyse the order conditions for an $s$-stages stiffly
accurate SRK method \eqref{SSRK}, let $\alpha=(\alpha_j) =
(A_{sj})^T$, let $\beta^{(k)}=(\beta^{(k)}_j)=(B^{(k)}_{sj})^T$ for
$k=1,2,3$ and define $e=(1,\ldots,1)^T \in \mathbb{R}^s$.
Because the stiffly accurate SRK method~\eqref{SSRK} is a special
case of the general class of SRK methods introduced in \cite{Roes},
the colored rooted tree theory in \cite{Roes,Roes2} can be applied
with Proposition~5.2 in \cite{Roes} to calculate order conditions
for the coefficients of the SRK method~\eqref{SSRK}. The strong
order 1.0 conditions for \eqref{SSRK} are calculated in \cite{KKR12}
and we print them here since we want to give a full classification
based on these order conditions in
Sections~\ref{Sec:Class-order-0.5} and \ref{Sec:Class-order-1.0}.
\begin{The} \label{SRK-theorem-ito-ord1.0-W1-main1}
    Let $f, g \in C^{1,3}(\mathcal{I} \times \mathbb{R}^d,
    \mathbb{R}^d)$.
    If the coefficients of the stochastic Runge-Kutta
    method~\eqref{SSRK}
    fulfill the equations
    {\allowdisplaybreaks
    \begin{alignat*}{3}
        1&. \quad {\alpha}^T e = 1
        \qquad \qquad \qquad \qquad \quad
        &2. \quad &{\beta^{(1)}}^T e = 1 \\
        3&. \quad {\beta^{(2)}}^T e = 0
        \qquad \qquad \qquad \qquad \quad
        &4. \quad &{\beta^{(3)}}^T e = 0 \\
        5&. \quad {\beta^{(1)}}^T B^{(1)} e = \frac{\lambda}{2}
        \qquad \qquad
        &6. \quad &{\beta^{(3)}}^T B^{(3)} e = -\frac{\lambda}{2} \\
        7&. \quad {\beta^{(2)}}^T B^{(3)} e + {\beta^{(3)}}^T B^{(2)} e =
      1-\lambda
        \qquad \qquad
        &8. \quad &{\alpha}^T B^{(3)} e = 0 \\
        9&. \quad {\beta^{(1)}}^T B^{(3)} e + {\beta^{(3)}}^T B^{(1)} e = 0
        \qquad \qquad
        &10. \quad &{\beta^{(2)}}^T B^{(2)} e = 0 \\
        11&. \quad {\beta^{(1)}}^T B^{(2)} e + {\beta^{(2)}}^T B^{(1)} e = 0
        \qquad \qquad
        &12. \quad &{\beta^{(3)}}^T A e = 0
    \end{alignat*}
    \begin{alignat*}{2}
        13. \quad &2 {\beta^{(1)}}^T (B^{(1)} e)(B^{(2)} e)
        + 2 {\beta^{(1)}}^T (B^{(1)} e)(B^{(3)} e)
        + {\beta^{(2)}}^T (B^{(1)} e)^2 \\
        & {} + {\beta^{(2)}}^T (B^{(2)} e)^2
        + {\beta^{(2)}}^T (B^{(2)} e)(B^{(3)} e)
        + {\beta^{(3)}}^T (B^{(1)} e)^2 \\
        & {} + \frac{1}{2} {\beta^{(3)}}^T (B^{(2)} e)^2
        + {\beta^{(3)}}^T (B^{(3)} e)^2 = 0 \\
        14. \quad &{\beta^{(1)}}^T (B^{(1)}(B^{(2)} e))
        + {\beta^{(1)}}^T (B^{(2)}(B^{(1)} e))
        + {\beta^{(1)}}^T (B^{(1)}(B^{(3)} e)) \\
        & {} + {\beta^{(1)}}^T (B^{(3)}(B^{(1)} e))
        + {\beta^{(2)}}^T (B^{(1)}(B^{(1)} e))
        + {\beta^{(2)}}^T (B^{(2)}(B^{(2)} e)) \\
        & {} + \frac{1}{2} {\beta^{(2)}}^T (B^{(2)}(B^{(3)} e))
        + \frac{1}{2} {\beta^{(2)}}^T (B^{(3)}(B^{(2)} e))
        + {\beta^{(3)}}^T (B^{(1)}(B^{(1)} e)) \\
        & {} + \frac{1}{2} {\beta^{(3)}}^T (B^{(2)}(B^{(2)} e))
        + {\beta^{(3)}}^T (B^{(3)}(B^{(3)} e)) = 0
    \end{alignat*}
    }
    for some $\lambda \in \mathbb{R}$ and
    if $c = A e$, then the stochastic Runge-Kutta
    method~\eqref{SSRK} attains order 1.0 for
    the strong approximation of the solution of the It{\^o}
    SDAE \eqref{SDAE-short} with scalar noise.
\end{The}
\begin{Bem} \label{Bem-Order-0.5-cond}
    Let $f, g \in C^{1,2}(\mathcal{I} \times \mathbb{R}^d,
    \mathbb{R}^d)$. Then, conditions 1--4 together with the
    condition ${\beta^{(1)}}^T B^{(1)} e + \frac{1}{2}
    {\beta^{(2)}}^T B^{(2)} e + {\beta^{(3)}}^T B^{(3)} e = 0$
    are sufficient for an order 0.5 strong SRK method \eqref{SSRK}
    that can be applied to the It{\^o} SDAE \eqref{SDAE-short}, see
    also \cite{KKR12}.
\end{Bem}
Using the order conditions, we will analyse the set of solutions in
the following
sections~\ref{Sec:Class-order-0.5}--\ref{Sec:Class-order-1.0}.
Because diagonally implicit SRK schemes are much more efficient with
respect to their computational effort compared to fully implicit SRK
schemes, we claim that $A_{ij}=B^{(3)}_{ij}=0$ for $j>i$ in the
following. Further, we need that $B^{(1)}_{ij}=B^{(2)}_{ij}=0$ for
$j \geq i$ in order to guarantee the existence of a solution for the
implicit equations in \eqref{SSRK} due to the unbounded random
variables $I_{(1),n}$ and $I_{(1,1),n}$, i.e., the SRK method has to
be explicit in the terms that involve random variables. Taking these
restrictions into account, we give a full classification for the
coefficients of the SRK method \eqref{SSRK}. Here, we have to point
out that in case of a singular matrix $M$ we choose the coefficients
within the classification such that either $A$ is regular or such
that $A_{1j}=B^{(k)}_{1j}=0$ and $A_{ii} \neq 0$ for $i\geq 2$.
Thus, the classification contains all coefficients such that the SRK
method \eqref{SSRK} can be applied to SDEs and may be explicit as
well as the case that it is implicit and can be applied to SDAEs.
Finally, the presented classification is the basis for the
calculation of coefficients for efficient SRK methods in the sense
that they primary possess a minimal number of stages, secondary have
a minimal number of implicit stages and finally for
Section~\ref{Sec:Optimal-Scheme-Stability-Ana} need a minimum of
explicit function evaluations. Under these restrictions, in
section~\ref{Sec:Optimal-Scheme-Stability-Ana} we try to find
efficient SRK schemes \eqref{SSRK} that are $A$-stable in the
mean-square sense.
\section{Classification of order 0.5 stiffly accurate SRK methods}
\label{Sec:Class-order-0.5}
Firstly, we give a full classification of strong order 0.5 stiffly
accurate SRK methods \eqref{SSRK} with a minimal number of stages
that can be diagonally implicit.
It easily follows that at least two stages are needed for
the order 0.5 conditions mentioned in Remark~\ref{Bem-Order-0.5-cond}
to be fulfilled. Therefore, $2$-stages SRK methods with coefficient table
\begin{equation} \label{sadirk_form}
    \begin{tabular}{|ccc|ccc|ccc|ccc}
        $A_{11}$ && &&& &&&& $B^{(3)}_{11}$ &&\\
        $A_{21}$ & $A_{22}$ & & $B^{(1)}_{21}$ & \hspace*{4mm} & & $B^{(2)}_{21}$ &
        \hspace*{4mm} & & $B^{(3)}_{21}$ & $B^{(3)}_{22}$ & \\
        \hline
    \end{tabular}
\end{equation}
are considered in this section. Because the considered SRK schemes
have to be explicit in terms involving random variables, the
coefficients $B^{(1)}_{11}$, $B^{(1)}_{22}$, $B^{(2)}_{11}$ and
$B^{(2)}_{22}$ are set equal to zero. Applying
Remark~\ref{Bem-Order-0.5-cond} to the case $s=2$ results in the
simplified system of order 0.5 conditions
\begin{align*}
1.\quad & A_{21}+A_{22} = 1 ,\\
2.\quad & B^{(1)}_{21} =1 ,\\
3.\quad & B^{(2)}_{21} =0 , \\
4.\quad & B^{(3)}_{21}+B^{(3)}_{22} =0 , \\
5.\quad & B^{(3)}_{21}B^{(3)}_{11} =0 .
\end{align*}
In the following, we denote by capital letters coefficients that can
be freely chosen whereas small letters stand for some prescribed
values. Solving these equations, we obviously get by simple
calculations the following two classes of order 0.5 stiffly accurate
SRK schemes \eqref{sadirk_form}:
\subsection{Strong order 0.5 SRK class I}
Choosing the coefficient $B^{(3)}_{21}=0$ implies that $B^{(3)}_{22}=0$ and defines class~I
with
\begin{equation} \label{Order-0.5-class-I}
    \begin{tabular}{|ccc|ccc|ccc|ccc}
     $A_{11}$ & & &\hspace*{4mm}&\hspace*{4mm}& \hspace*{4mm}
     &\hspace*{4mm}&\hspace*{4mm}&\hspace*{4mm}& $B_{11}^{(3)}$ &&\\
     $A_{21}$ & $a_{22}$ & & $1$ & & & $0$ & && $0$
     & $0$ & \\
     \hline
   \end{tabular}
\end{equation}
where $a_{22}=1-A_{21}$ and $A_{11}, A_{21}, B_{11}^{(3)} \in
\mathbb{R}$. \\ \\
Remark that in case of $A_{11}=B_{11}^{(3)}=0$, the SRK scheme
\eqref{SSRK} with coefficients \eqref{Order-0.5-class-I} coincides
with the well known stochastic $\theta$-method in \cite{Hig00}.
\subsection{Strong order 0.5 SRK class II}
If we choose $B_{11}^{(3)}=0$, then we get the class~II coefficients with
\begin{equation} \label{Order-0.5-class-II}
    \begin{tabular}{|ccc|ccc|ccc|ccc}
     $A_{11}$ & & &\hspace*{4mm}&\hspace*{4mm}& \hspace*{4mm}
     &\hspace*{4mm}&\hspace*{4mm}&\hspace*{4mm}& $0$ &&\\
     $A_{21}$ & $a_{22}$ & & $1$ & & & $0$ & && $B_{21}^{(3)}$
     & $b_{22}^{(3)}$ & \\
     \hline
   \end{tabular}
\end{equation}
where $a_{22}=1-A_{21}$, $b_{22}^{(3)}=-B_{21}^{(3)}$ and $A_{11},
A_{21}, B_{21}^{(3)} \in \mathbb{R}$.
\section{Classification of order 1.0 stiffly accurate SRK methods}
\label{Sec:Class-order-1.0}
Next, we search for stiffly accurate diagonally implicit SRK methods of
strong order 1.0 with a minimal number of stages. Again, these
methods should be explicit in the terms involving random variables.
From the order 1.0 conditions given in Theorem~\ref{SRK-theorem-ito-ord1.0-W1-main1}
it follows that a minimum number of $s=3$ stages are required. This can be seen
easily, because for some smaller $s$ there exist no coefficients that fulfill
the order conditions 2, 3, 5 and 7 in Theorem~\ref{SRK-theorem-ito-ord1.0-W1-main1}.
Thus, at least $s=3$ stages are needed to assure strong order 1.0 for the SRK method.
These $3$-stages stiffly accurate diagonally implicit SRK schemes are determined
by the following coefficient table:
\begin{equation} \label{sadirk_form1.0}
\begin{tabular}{|ccc|ccc|ccc|ccc}
  $A_{11}$ && &&& &&&& $B^{(3)}_{11}$ &&\\
  $A_{21}$ & $A_{22}$ & & $B^{(1)}_{21}$ & & & $B^{(2)}_{21}$ & &&
  $B^{(3)}_{21}$ & $B^{(3)}_{22}$ & \\
  $A_{31}$ & $A_{32}$ & $A_{33}$ & $B^{(1)}_{31}$ & $B^{(1)}_{32}$ & \hspace*{4mm}&
  $B^{(2)}_{31}$ & $B^{(2)}_{32}$ & \hspace*{4mm} & $B^{(3)}_{31}$ & $B^{(3)}_{32}$ &
  $B^{(3)}_{33}$ \\ \hline
\end{tabular}
\end{equation}
Then, the first four order conditions of Theorem~\ref{SRK-theorem-ito-ord1.0-W1-main1}
reduce to
\begin{align*}
1.\quad & A_{31}+A_{32}+A_{33} = 1 , \\
2.\quad & B^{(1)}_{31}+B^{(1)}_{32} =1 , \\
3.\quad & B^{(2)}_{31}+B^{(2)}_{32} =0 , \\
4.\quad & B^{(3)}_{31}+B^{(3)}_{32}+B^{(3)}_{33} =0 .
\end{align*}
Taking into account these simplified conditions, the remaining conditions 5--12 can
be written as
\begin{align*}
5.\quad & B^{(1)}_{32}B^{(1)}_{21} = \frac{\lambda}{2} , \\
6.\quad & B^{(3)}_{31}B^{(3)}_{11} + B^{(3)}_{32}(B^{(3)}_{21} + B^{(3)}_{22}) = -\frac{\lambda}{2} , \\
7.\quad & B^{(2)}_{31}B^{(3)}_{11} + B^{(2)}_{32}(B^{(3)}_{21} + B^{(3)}_{22}) + B^{(3)}_{32}B^{(2)}_{21}  = 1-\lambda , \\
8.\quad & A_{31}B^{(3)}_{11} + A_{32}(B^{(3)}_{21} + B^{(3)}_{22}) = 0 , \\
9.\quad & B^{(1)}_{31}B^{(3)}_{11} + B^{(1)}_{32}(B^{(3)}_{21} + B^{(3)}_{22}) + B^{(3)}_{32}B^{(1)}_{21} + B^{(3)}_{33} = 0 , \\
10.\quad & B^{(2)}_{32}B^{(2)}_{21} = 0 , \\
11.\quad & B^{(1)}_{32}B^{(2)}_{21} + B^{(2)}_{32}B^{(1)}_{21} = 0 , \\
12.\quad & B^{(3)}_{31}A_{11} + B^{(3)}_{32}(A_{21}+A_{22}) + B^{(3)}_{33} = 0 .
\end{align*}
For conditions 13 and 14 we refer to Theorem~\ref{SRK-theorem-ito-ord1.0-W1-main1}. Then,
the following result can be derived in the case of $s=3$ from the simplified order conditions.
\begin{Lem}
For a stiffly accurate order 1.0 SRK method \eqref{SSRK} with three stages and coefficient
scheme \eqref{sadirk_form1.0} the following assertions hold:
\begin{enumerate}[(i)]
\item For the parameter $\lambda$ which occurs in the order conditions 5--7 follows that
  $\lambda \in \{0,1\}$.
\item It holds $\lambda=1$ if and only if $B^{(2)}=0$.
\end{enumerate}
\end{Lem}
{\bf Proof.}
The results follow straight forward from the solution of the order conditions: Assume that
$\lambda \neq 0$. Then condition 5 yields, that $B^{(1)}_{32}\neq 0$ and $B^{(1)}_{21}\neq 0$.
From condition 10 we get, that $B^{(2)}_{32}=0$ or $B^{(2)}_{21}=0$ and therefore at least
one of the terms on the left hand side of condition 11 is equal to 0. Then the other term
on the left hand side of condition 11 also has to be 0 and thus $B^{(2)}_{32}=0$ and
$B^{(2)}_{21}=0$. Condition 3 yields, that $B^{(2)}_{31}=0$. Therefore we have $B^{(2)}=0$.
Now, the left hand side of condition 7 vanishes, thus we get $\lambda=1$. This proves (i)
and (ii). \hfill $\Box$ \\ \\
For the analysis of the set of coefficients that fulfill the strong
order 1.0 conditions, we derive the following possible classes of
schemes, where we have $\lambda=1$ for the first five classes and
$\lambda=0$ for the remaining six classes. Most of the calculations
are done using the software Maple. All presented classes are
significantly different although not totally disjoint due to our
choice of a clear and compact way for their representation. Special
attention has to be paid to the signs of some of the coefficients.
Whenever positive as well as negative signs are allowed, one has to
choose either the upper or the lower sign of the symbols $\pm$ and
$\mp$, respectively, for all affected coefficients.
In the following, we denote all coefficients that can be chosen
freely by capital letters, whereas lower case is used to denote more
complex expressions.
\subsection{Strong order 1.0 SRK class I with $\lambda=1$}
The first class of coefficients is given for $A_{11}, A_{22},
A_{33}, B_{22}^{(3)} \in \mathbb{R}$ and $B_{32}^{(3)} \in
\mathbb{R} \setminus \{0\}$ by the tableau
\begin{equation}
    \begin{tabular}{|ccc|ccc|ccc|ccc}
     $A_{11}$ & & &\hspace*{5mm}&\hspace*{5mm}& \hspace*{5mm}
     &\hspace*{5mm}&\hspace*{5mm}&\hspace*{5mm}& $0$ &&\\
     $a_{21}$ & $A_{22}$ & & $1$ & & & $0$ & && $b_{21}^{(3)}$
     & $B_{22}^{(3)}$ & \\
     $a_{31}$ & $0$ & $A_{33}$ & $\frac{1}{2}$ & $\frac{1}{2}$& &
     $0$ & $0$ & & $b_{31}^{(3)}$ & $B_{32}^{(3)}$ &
     $b_{33}^{(3)}$ \\ \hline
   \end{tabular}
\end{equation}
with fixed coefficients
\begin{alignat*}{3}
    a_{21} &= \tfrac{A_{11}-4 A_{22} (B_{32}^{(3)})^2 + 4 (B_{32}^{(3)})^2
    -1}{4 (B_{32}^{(3)})^2},
    &\quad \quad
    a_{31} &= 1-A_{33},
    &\quad \quad
    b_{21}^{(3)} &= -\tfrac{1 + 2 B_{22}^{(3)} B_{32}^{(3)}}{2
    B_{32}^{(3)}}, \\
    b_{31}^{(3)} &= -\tfrac{1}{4 B_{32}^{(3)}},
    &\quad \quad
    b_{33}^{(3)} &= -\tfrac{4 (B_{32}^{(3)})^2 -1}{4 B_{32}^{(3)}}.
\end{alignat*}
\subsection{Strong order 1.0 SRK class II with $\lambda=1$}
The second class is given for $A_{11}, A_{22}, A_{33}, B_{22}^{(3)}
\in \mathbb{R}$ and $B_{32}^{(3)} \in \mathbb{R} \setminus \{0\}$ by
the tableau
\begin{equation}
    \begin{tabular}{|ccc|ccc|ccc|ccc}
     $A_{11}$ & & &\hspace*{5mm}&\hspace*{5mm}& \hspace*{5mm}
     &\hspace*{5mm}&\hspace*{5mm}&\hspace*{5mm}& $0$ &&\\
     $a_{21}$ & $A_{22}$ & & $b_{21}^{(1)}$ & & & $0$ & && $b_{21}^{(3)}$
     & $B_{22}^{(3)}$ & \\
     $a_{31}$ & $0$ & $A_{33}$ & $b_{31}^{(1)}$ & $b_{32}^{(1)}$& &
     $0$ & $0$ & & $b_{31}^{(3)}$ & $B_{32}^{(3)}$ &
     $0$ \\ \hline
   \end{tabular}
\end{equation}
with fixed coefficients
\begin{alignat*}{4}
    a_{21} &= A_{11}-A_{22}, &\quad \quad
    a_{31} &= 1-A_{33}, &\quad \quad
    b_{21}^{(1)} &= \pm \tfrac{1}{2 B_{32}^{(3)}}, &\quad \quad
    b_{31}^{(1)} &= 1 \mp B_{32}^{(3)}, \\
    b_{32}^{(1)} &= \pm B_{32}^{(3)}, &\quad \quad
    b_{21}^{(3)} &= -\tfrac{1+2 B_{22}^{(3)} B_{32}^{(3)}}{2
    B_{32}^{(3)}}, &\quad \quad
    b_{31}^{(3)} &= -B_{32}^{(3)}.
\end{alignat*}
\subsection{Strong order 1.0 SRK class III with $\lambda=1$}
The third class of coefficients is determined for $A_{21}, A_{22},
A_{32} \in \mathbb{R}$ and $B_{11}^{(3)} \in \mathbb{R} \setminus
\{0\}$ by the tableau
\begin{equation}
    \begin{tabular}{|ccc|ccc|ccc|ccc}
     $1$ & & &\hspace*{5mm}&\hspace*{5mm}& \hspace*{5mm}
     &\hspace*{5mm}&\hspace*{5mm}&\hspace*{5mm}& $B_{11}^{(3)}$ &&\\
     $A_{21}$ & $A_{22}$ & & $b_{21}^{(1)}$ & & & $0$ & && $b_{21}^{(3)}$
     & $0$ & \\
     $a_{31}$ & $A_{32}$ & $a_{33}$ & $b_{31}^{(1)}$ & $b_{32}^{(1)}$& &
     $0$ & $0$ & & $b_{31}^{(3)}$ & $0$ &
     $b_{33}^{(3)}$ \\ \hline
   \end{tabular}
\end{equation}
with fixed coefficients
\begin{alignat*}{4}
    a_{31} &= -\tfrac{A_{32} ((B_{11}^{(3)})^2 -1)}{2 (B_{11}^{(3)})^2},
    &\quad \quad
    a_{33} &= -\tfrac{A_{32} (B_{11}^{(3)})^2 -2 (B_{11}^{(3)})^2 +
    A_{32}}{2 (B_{11}^{(3)})^2},
    &\quad \quad
    b_{21}^{(1)} &= \tfrac{1}{2} \tfrac{(B_{11}^{(3)})^2+1}{1+2 (B_{11}^{(3)})^2}, \\
    b_{31}^{(1)} &= -\tfrac{(B_{11}^{(3)})^2}{(B_{11}^{(3)})^2 +1},
    &\quad \quad
    b_{32}^{(1)} &= \tfrac{1+ 2(B_{11}^{(3)})^2}{(B_{11}^{(3)})^2 +1},
    &\quad \quad
    b_{21}^{(3)} &= \tfrac{(B_{11}^{(3)})^2 -1}{2 B_{11}^{(3)}}, \\
    b_{31}^{(3)} &= -\tfrac{1}{2 B_{11}^{(3)}},
    &\quad \quad
    b_{33}^{(3)} &= \tfrac{1}{2 B_{11}^{(3)}}.
\end{alignat*}
\subsection{Strong order 1.0 SRK class IV with $\lambda=1$}
For the fourth class, for $A_{11}, A_{22}, A_{32} \in \mathbb{R}$
and $B_{33}^{(3)} \in \mathbb{R} \setminus \{0\}$ the coefficients
are given by the tableau
\begin{equation}
    \begin{tabular}{|ccc|ccc|ccc|ccc}
     $A_{11}$ & & &\hspace*{5mm}&\hspace*{5mm}& \hspace*{5mm}
     &\hspace*{5mm}&\hspace*{5mm}&\hspace*{5mm}& $b_{11}^{(3)}$ &&\\
     $a_{21}$ & $A_{22}$ & & $b_{21}^{(1)}$ & & & $0$ & && $b_{21}^{(3)}$
     & $0$ & \\
     $a_{31}$ & $A_{32}$ & $a_{33}$ & $b_{31}^{(1)}$ & $b_{32}^{(1)}$& &
     $0$ & $0$ & & $b_{31}^{(3)}$ & $b_{32}^{(3)}$ &
     $B_{33}^{(3)}$ \\ \hline
   \end{tabular}
\end{equation}
with fixed coefficients
\begin{alignat*}{4}
    a_{21} &= \tfrac{2 (B_{33}^{(3)})^2 -2 A_{22} (B_{33}^{(3)})^2 +2
    -A_{22} -A_{11}}{1+2 (B_{33}^{(3)})^2},
    \ \ 
    a_{31} = \tfrac{A_{32}}{(B_{33}^{(3)})^2},
    \ \ 
    a_{33} = \tfrac{(B_{33}^{(3)})^2 -(B_{33}^{(3)})^2 A_{32}
    -A_{32}}{(B_{33}^{(3)})^2}, \\
    b_{21}^{(1)} &= \pm \tfrac{1+(B_{33}^{(3)})^2}{B_{33}^{(3)}
    \sqrt{1+2 (B_{33}^{(3)})^2}},
    \quad 
    b_{31}^{(1)} = -\tfrac{1}{2} \tfrac{\pm B_{33}^{(3)}
    \sqrt{1+2(B_{33}^{(3)})^2}-2-2
    (B_{33}^{(3)})^2}{1+(B_{33}^{(3)})^2}, \\
    b_{32}^{(1)} &= \pm \tfrac{1}{2} \tfrac{B_{33}^{(3)}
    \sqrt{1+2 (B_{33}^{(3)})^2}}{1+ (B_{33}^{(3)})^2},
    \quad 
    b_{11}^{(3)} = -B_{33}^{(3)},
    \quad 
    b_{21}^{(3)} = \tfrac{1}{B_{33}^{(3)}},
    \quad 
    b_{31}^{(3)} = -\tfrac{1}{2}
    \tfrac{B_{33}^{(3)}}{1+(B_{33}^{(3)})^2}, \\
    b_{32}^{(3)} &= -\tfrac{1}{2} \tfrac{B_{33}^{(3)}
    (1+2 (B_{33}^{(3)})^2)}{1+(B_{33}^{(3)})^2}.
\end{alignat*}
\subsection{Strong order 1.0 SRK class V with $\lambda=1$}
The fifth class of coefficients is defined for $A_{11}, A_{22},
A_{32} \in \mathbb{R}$ by the tableau
\begin{equation}
    \begin{tabular}{|ccc|ccc|ccc|ccc}
     $A_{11}$ & & &\hspace*{5mm}&\hspace*{5mm}& \hspace*{5mm}
     &\hspace*{5mm}&\hspace*{5mm}&\hspace*{5mm}& $b_{11}^{(3)}$ &&\\
     $a_{21}$ & $A_{22}$ & & $b_{21}^{(1)}$ & & & $0$ & && $b_{21}^{(3)}$
     & $0$ & \\
     $a_{31}$ & $A_{32}$ & $a_{33}$ & $b_{31}^{(1)}$ & $B_{32}^{(1)}$& &
     $0$ & $0$ & & $b_{31}^{(3)}$ & $B_{32}^{(3)}$ &
     $B_{33}^{(3)}$ \\ \hline
   \end{tabular}
\end{equation}
with fixed coefficients
\begin{alignat*}{4}
    a_{21} &= \tfrac{-B_{33}^{(3)} +A_{11} B_{33}^{(3)} -A_{22}
    B_{32}^{(3)} +A_{11} B_{32}^{(3)}}{B_{32}^{(3)}},
    \quad \\
    a_{31} &= -\tfrac{((B_{32}^{(1)})^2 -2 B_{32}^{(1)} B_{33}^{(3)}
    B_{32}^{(3)} -2B_{32}^{(1)} (B_{33}^{(3)})^2 -B_{32}^{(1)}
    -B_{32}^{(3)} B_{33}^{(3)} -(B_{32}^{(3)})^2)
    A_{32}}{(B_{32}^{(1)})^2 -2B_{32}^{(1)} B_{33}^{(3)} B_{32}^{(3)}
    -(B_{32}^{(3)})^2},
    \quad \\
    a_{33} &= -\tfrac{-(B_{32}^{(1)})^2 +2B_{32}^{(1)} B_{33}^{(3)}
    B_{32}^{(3)} +2B_{32}^{(1)} A_{32} (B_{33}^{(3)})^2 +A_{32}
    B_{32}^{(1)} +A_{32} B_{33}^{(3)} B_{32}^{(3)}
    +(B_{32}^{(3)})^2}{(B_{32}^{(1)})^2 -2B_{32}^{(1)} B_{33}^{(3)}
    B_{32}^{(3)} -(B_{32}^{(3)})^2},
    \quad \\
    b_{21}^{(1)} &= \tfrac{1}{2 B_{32}^{(1)}},
    \quad
    b_{31}^{(1)}=1-B_{32}^{(1)},
    \quad
    b_{11}^{(3)}=\tfrac{1}{2} \tfrac{(B_{32}^{(1)})^2 -2B_{32}^{(1)}
    B_{33}^{(3)} B_{32}^{(3)} -(B_{32}^{(3)})^2}{B_{32}^{(1)}
    (B_{32}^{(1)} B_{33}^{(3)} +B_{32}^{(3)})}, \\
    b_{21}^{(3)} &= \tfrac{1}{2} \tfrac{(B_{32}^{(1)})^2 -2B_{32}^{(1)}
    B_{33}^{(3)} B_{32}^{(3)} -2B_{32}^{(1)} (B_{33}^{(3)})^2
    -B_{32}^{(1)} -B_{32}^{(3)} B_{33}^{(3)}
    -(B_{32}^{(3)})^2}{B_{32}^{(1)} (B_{32}^{(1)} B_{33}^{(3)}
    +B_{32}^{(3)})},
    \ \ 
    b_{31}^{(3)}=-B_{32}^{(3)}-B_{33}^{(3)}
\end{alignat*}
and all solutions $B_{32}^{(1)}, B_{32}^{(3)}, B_{33}^{(3)} \in \mathbb{R}
\setminus \{0\}$ of the equation
\begin{equation} \label{Order-cond-implicit-eqn-lambda1}
    \begin{split}
    &4 (B_{32}^{(1)})^2 B_{32}^{(3)} (B_{33}^{(3)})^4
    + 4 B_{32}^{(1)} (B_{32}^{(3)})^2 (B_{33}^{(3)})^3
    + 4 (B_{32}^{(1)})^3 B_{32}^{(3)} (B_{33}^{(3)})^2 \\
    &+ 4 (B_{32}^{(1)})^2 (B_{32}^{(3)})^2 (B_{33}^{(3)})^3
    + 4 B_{32}^{(1)} (B_{32}^{(3)})^3 (B_{33}^{(3)})^2
    - 4 (B_{32}^{(1)})^3 (B_{33}^{(3)})^3 \\
    &+ (B_{32}^{(3)})^3 (B_{33}^{(3)})^2
    - 2 (B_{32}^{(1)})^3 B_{33}^{(3)}
    - (B_{32}^{(1)})^2 B_{32}^{(3)}
    - (B_{32}^{(1)})^2 B_{32}^{(3)} (B_{33}^{(3)})^2 \\
    &+ (B_{32}^{(3)})^3
    + 2 (B_{32}^{(1)})^2 (B_{32}^{(3)})^2 B_{33}^{(3)}
    + 2 B_{32}^{(1)} (B_{32}^{(3)})^2 B_{33}^{(3)}
    + (B_{32}^{(1)})^4 B_{33}^{(3)} \\
    &+ (B_{32}^{(3)})^4 B_{33}^{(3)}
    + 4 (B_{32}^{(1)})^4 (B_{33}^{(3)})^3 = 0
    \end{split}
\end{equation}
where 
$B_{32}^{(3)} \neq -B_{32}^{(1)} B_{33}^{(3)}$ is needed.
\subsection{Strong order 1.0 SRK class VI with $\lambda=0$}
For $\lambda=0$, class six is given by the coefficients $A_{11},
A_{22}, A_{32} \in \mathbb{R}$ with the tableau
\begin{equation}
    \begin{tabular}{|ccc|ccc|ccc|ccc}
     $A_{11}$ & & &\hspace*{5mm}&\hspace*{5mm}& \hspace*{5mm}
     &\hspace*{5mm}&\hspace*{5mm}&\hspace*{5mm}& $B_{11}^{(3)}$ &&\\
     $a_{21}$ & $A_{22}$ & & $b_{21}^{(1)}$ & & & $b_{21}^{(2)}$ & && $b_{21}^{(3)}$
     & $0$ & \\
     $a_{31}$ & $A_{32}$ & $a_{33}$ & $1$ & $0$& &
     $0$ & $0$ & & $b_{31}^{(3)}$ & $B_{32}^{(3)}$ &
     $b_{33}^{(3)}$ \\ \hline
   \end{tabular}
\end{equation}
where
\begin{alignat*}{4}
    a_{21} &= -\tfrac{1}{2} \tfrac{1}{((B_{11}^{(3)})^2 +1) B_{32}^{(3)}}
    (-A_{11} B_{32}^{(3)} (B_{11}^{(3)})^2 -A_{11} B_{32}^{(3)} +2
    A_{11} B_{11}^{(3)} \pm A_{11} \sqrt{D} \\
    & {} \quad +2A_{22} B_{32}^{(3)}
    (B_{11}^{(3)})^2 +2A_{22} B_{32}^{(3)} -2B_{11}^{(3)}
    -(B_{11}^{(3)})^2 B_{32}^{(3)} -B_{32}^{(3)} \mp \sqrt{D}), \\
    a_{31} &= \tfrac{1}{2} \tfrac{1}{((B_{11}^{(3)})^2 +1) B_{32}^{(3)}}
    (A_{32} (-(B_{11}^{(3)})^2 B_{32}^{(3)} -B_{32}^{(3)} +2B_{11}^{(3)}
    \pm \sqrt{D})), \\
    a_{33} &= -\tfrac{1}{2} \tfrac{1}{((B_{11}^{(3)})^2 +1)
    B_{32}^{(3)}} (-2(B_{11}^{(3)})^2 B_{32}^{(3)} -2B_{32}^{(3)}
    +2B_{11}^{(3)} A_{32} +A_{32} (B_{11}^{(3)})^2 B_{32}^{(3)} \\
    & {} \quad +A_{32} B_{32}^{(3)} \pm A_{32} \sqrt{D}), \\
    b_{21}^{(1)} &= \tfrac{1}{2} \tfrac{1}{((B_{11}^{(3)})^2 +1)
    B_{32}^{(3)}} ((B_{11}^{(3)})^2
    B_{32}^{(3)} +B_{32}^{(3)} -2(B_{11}^{(3)})^3 \pm \sqrt{D}),
    \quad
    b_{21}^{(2)} = \tfrac{1}{B_{32}^{(3)}}, \\
    b_{21}^{(3)} &= -\tfrac{1}{2} \tfrac{1}{((B_{11}^{(3)})^2 +1)
    B_{32}^{(3)}} (B_{11}^{(3)} (-(B_{11}^{(3)})^2 B_{32}^{(3)}
    -B_{32}^{(3)} +2B_{11}^{(3)} \pm \sqrt{D})), \\
    b_{31}^{(3)} &= \tfrac{1}{2} \tfrac{1}{(B_{11}^{(3)})^2 +1}
    (-(B_{11}^{(3)})^2 B_{32}^{(3)} -B_{32}^{(3)} +2B_{11}^{(3)} \pm
    \sqrt{D}), \\
    b_{33}^{(3)} &= -\tfrac{1}{2} \tfrac{1}{(B_{11}^{(3)})^2 +1}
    (2B_{11}^{(3)} +(B_{11}^{(3)})^2 B_{32}^{(3)} +B_{32}^{(3)} \pm
    \sqrt{D})
\end{alignat*}
with $B_{11}^{(3)} \in \mathbb{R}$, $B_{32}^{(3)} \in \mathbb{R}
\setminus \{0\}$ and
\begin{equation}
    \begin{split}
    D &= (B_{11}^{(3)})^4 (B_{32}^{(3)})^2 +2 (B_{11}^{(3)})^2
    (B_{32}^{(3)})^2 +(B_{32}^{(3)})^2 +4(B_{11}^{(3)})^3
    B_{32}^{(3)} \\
    & {} \quad -2(B_{11}^{(3)})^2
    -4(B_{11}^{(3)})^4 +4B_{11}^{(3)} B_{32}^{(3)} -2
    \end{split}
\end{equation}
such that $D \geq 0$ is fulfilled.
\subsection{Strong order 1.0 SRK class VII with $\lambda=0$}
Class seven is defined for $A_{11}, A_{22}, A_{32}, A_{33},
B_{22}^{(3)} \in \mathbb{R}$ and $B_{21}^{(1)} \in \mathbb{R}
\setminus \{0\}$ by the tableau
\begin{equation}
    \begin{tabular}{|ccc|ccc|ccc|ccc}
     $A_{11}$ & & &\hspace*{5mm}&\hspace*{5mm}& \hspace*{5mm}
     &\hspace*{5mm}&\hspace*{5mm}&\hspace*{5mm}& $0$ &&\\
     $a_{21}$ & $A_{22}$ & & $B_{21}^{(1)}$ & & & $b_{21}^{(2)}$ & && $b_{21}^{(3)}$
     & $B_{22}^{(3)}$ & \\
     $a_{31}$ & $A_{32}$ & $A_{33}$ & $1$ & $0$& &
     $0$ & $0$ & & $b_{31}^{(3)}$ & $b_{32}^{(3)}$ &
     $b_{33}^{(3)}$ \\ \hline
   \end{tabular}
\end{equation}
where
\begin{alignat*}{4}
    a_{21} &= A_{11} -A_{11} B_{21}^{(1)} -A_{22} +B_{21}^{(1)},
    &\quad \quad \quad
    a_{31} &= 1-A_{32}-A_{33}, \\
    b_{21}^{(2)} &= \pm \sqrt{2 B_{21}^{(1)} -2 (B_{21}^{(1)})^2},
    &\quad \quad
    b_{21}^{(3)} &= -B_{22}^{(3)}, \\
    b_{31}^{(3)} &= \mp \frac{ 1 -B_{21}^{(1)}}{\sqrt{2B_{21}^{(1)}
    -2(B_{21}^{(1)})^2}},
    &\quad \quad
    b_{32}^{(3)} &= \pm
    \frac{1}{\sqrt{2B_{21}^{(1)}-2(B_{21}^{(1)})^2}}, \\
    b_{33}^{(3)} &= \mp \frac{B_{21}^{(1)}}{\sqrt{2B_{21}^{(1)}
    -2(B_{21}^{(1)})^2}}.
\end{alignat*}
\subsection{Strong order 1.0 SRK class VIII with $\lambda=0$}
For $A_{11}, A_{21}, A_{22}, A_{32}, B_{22}^{(3)} \in \mathbb{R}$
and $B_{32}^{(2)}, B_{11}^{(3)} \in \mathbb{R} \setminus \{0\}$, the
eighth class is given by the tableau
\begin{equation}
    \begin{tabular}{|ccc|ccc|ccc|ccc}
     $A_{11}$ & & &\hspace*{5mm}&\hspace*{5mm}& \hspace*{5mm}
     &\hspace*{5mm}&\hspace*{5mm}&\hspace*{5mm}& $B_{11}^{(3)}$ &&\\
     $A_{21}$ & $A_{22}$ & & $0$ & & & $0$ & && $b_{21}^{(3)}$
     & $B_{22}^{(3)}$ & \\
     $a_{31}$ & $A_{32}$ & $a_{33}$ & $b_{31}^{(1)}$ & $b_{32}^{(1)}$& &
     $b_{31}^{(2)}$ & $B_{32}^{(2)}$ & & $0$ & $0$ &
     $0$ \\ \hline
   \end{tabular}
\end{equation}
with
\begin{alignat*}{4}
    a_{31} &= -\frac{A_{32} (1+B_{32}^{(2)} B_{11}^{(3)})}{B_{32}^{(2)}
    B_{11}^{(3)}},
    \quad 
    a_{33} = \frac{B_{32}^{(2)} B_{11}^{(3)} +A_{32}}{B_{32}^{(2)}
    B_{11}^{(3)}},
    \quad 
    b_{31}^{(1)} = 1 +B_{32}^{(2)} B_{11}^{(3)}, \\
    b_{32}^{(1)} &= -B_{32}^{(2)} B_{11}^{(3)},
    \quad 
    b_{31}^{(2)} = -B_{32}^{(2)},
    \quad 
    b_{21}^{(3)} = \frac{1 +B_{32}^{(2)} (B_{11}^{(3)}
    -B_{22}^{(3)})}{B_{32}^{(2)}}.
\end{alignat*}
\subsection{Strong order 1.0 SRK class IX with $\lambda=0$}
Class nine with $\lambda=0$ is given for $A_{11}, A_{22}, A_{32},
B_{32}^{(3)} \in \mathbb{R}$ and $B_{11}^{(3)} \in \mathbb{R}
\setminus \{0\}$ by the tableau
\begin{equation}
    \begin{tabular}{|ccc|ccc|ccc|ccc}
     $A_{11}$ & & &\hspace*{5mm}&\hspace*{5mm}& \hspace*{5mm}
     &\hspace*{5mm}&\hspace*{5mm}&\hspace*{5mm}& $B_{11}^{(3)}$ &&\\
     $a_{21}$ & $A_{22}$ & & $0$ & & & $0$ & && $b_{21}^{(3)}$
     & $0$ & \\
     $a_{31}$ & $A_{32}$ & $a_{33}$ & $b_{31}^{(1)}$ & $b_{32}^{(1)}$& &
     $b_{31}^{(2)}$ & $b_{32}^{(2)}$ & & $b_{31}^{(3)}$ & $B_{32}^{(3)}$ &
     $b_{33}^{(3)}$ \\ \hline
   \end{tabular}
\end{equation}
with the coefficients
\begin{alignat*}{4}
    a_{21} &= \frac{(B_{11}^{(3)})^2 -A_{22} (B_{11}^{(3)})^2 -A_{11}
    +1}{(B_{11}^{(3)})^2},
    &\quad \quad
    a_{31} &= \frac{A_{32}}{(B_{11}^{(3)})^2}, \\
    a_{33} &= \frac{(B_{11}^{(3)})^2 -(B_{11}^{(3)})^2 A_{32}
    -A_{32}}{(B_{11}^{(3)})^2},
    &\quad
    b_{31}^{(1)} &= \frac{B_{11}^{(3)}
    +(B_{11}^{(3)})^2 B_{32}^{(3)} +B_{32}^{(3)}}{B_{11}^{(3)}
    ((B_{11}^{(3)})^2 +1)}, \\
    b_{32}^{(1)} &= \frac{(B_{11}^{(3)})^3
    -(B_{11}^{(3)})^2 B_{32}^{(3)} -B_{32}^{(3)}}{B_{11}^{(3)}
    ((B_{11}^{(3)})^2 +1)},
    &\quad
    b_{31}^{(2)} &= \frac{B_{11}^{(3)}}{(B_{11}^{(3)})^2 +1}, \\
    b_{32}^{(2)} &= -\frac{B_{11}^{(3)}}{(B_{11}^{(3)})^2 +1},
    &\quad
    b_{21}^{(3)} &= -\frac{1}{B_{11}^{(3)}}, \\
    b_{31}^{(3)} &= \frac{B_{32}^{(3)}}{(B_{11}^{(3)})^2},
    &\quad
    b_{33}^{(3)} &= -\frac{((B_{11}^{(3)})^2 +1)
    B_{32}^{(3)}}{(B_{11}^{(3)})^2}.
\end{alignat*}
\subsection{Strong order 1.0 SRK class X with $\lambda=0$}
Class ten is defined for $A_{11}, A_{21}, A_{22}, A_{33}, B_{22}^{(3)} \in
\mathbb{R}$ and $B_{32}^{(2)} \in \mathbb{R} \setminus \{0\}$ by the tableau
\begin{equation}
    \begin{tabular}{|ccc|ccc|ccc|ccc}
     $A_{11}$ & & &\hspace*{5mm}&\hspace*{5mm}& \hspace*{5mm}
     &\hspace*{5mm}&\hspace*{5mm}&\hspace*{5mm}& $0$ &&\\
     $A_{21}$ & $A_{22}$ & & $0$ & & & $0$ & && $b_{21}^{(3)}$
     & $B_{22}^{(3)}$ & \\
     $a_{31}$ & $0$ & $A_{33}$ & $1$ & $0$& &
     $b_{31}^{(2)}$ & $B_{32}^{(2)}$ & & $0$ & $0$ &
     $0$ \\ \hline
   \end{tabular}
\end{equation}
with coefficients
\begin{alignat*}{4}
    a_{31} &= 1-A_{33},
    \quad \quad
    b_{31}^{(2)} = -B_{32}^{(2)},
    \quad \quad
    b_{21}^{(3)} = \frac{1 -B_{32}^{(2)} B_{22}^{(3)}}{B_{32}^{(2)}} .
\end{alignat*}
\subsection{Strong order 1.0 SRK class XI with $\lambda=0$}
The last class eleven is given for $A_{21}, A_{22}, A_{33} \in \mathbb{R}$ and
$B_{33}^{(3)} \in \mathbb{R} \setminus \{0\}$ by the tableau
\begin{equation}
    \begin{tabular}{|ccc|ccc|ccc|ccc}
     $a_{11}$ & & &\hspace*{5mm}&\hspace*{5mm}& \hspace*{5mm}
     &\hspace*{5mm}&\hspace*{5mm}&\hspace*{5mm}& $b_{11}^{(3)}$ &&\\
     $A_{21}$ & $A_{22}$ & & $0$ & & & $b_{21}^{(2)}$ & && $b_{21}^{(3)}$
     & $0$ & \\
     $a_{31}$ & $a_{32}$ & $A_{33}$ & $1$ & $0$& &
     $0$ & $0$ & & $b_{31}^{(3)}$ & $b_{32}^{(3)}$ &
     $B_{33}^{(3)}$ \\ \hline
   \end{tabular}
\end{equation}
and the coefficients
\begin{alignat*}{4}
  a_{11} &= \frac{2A_{21} (B_{33}^{(3)})^4 +A_{21} +2A_{22}
  (B_{33}^{(3)})^4 +A_{22} -2(B_{33}^{(3)})^4
  -2(B_{33}^{(3)})^2}{1-2(B_{33}^{(3)})^2}, \\
  a_{31} &= 1-a_{32}-A_{33}, \quad \quad \quad
  a_{32} = -\frac{2 A_{33} (B_{33}^{(3)})^4
  -2(B_{33}^{(3)})^4 -1 +A_{33}}{2 (B_{33}^{(3)})^2 ((B_{33}^{(3)})^2
  +1)}, \\
  b_{21}^{(2)} &= \frac{\sqrt{-2(B_{33}^{(3)})^3 b_{32}^{(3)}
  -2b_{32}^{(3)} B_{33}^{(3)} -2(B_{33}^{(3)})^4}}{b_{32}^{(3)}}, \\
  b_{11}^{(3)} &= -B_{33}^{(3)}, \quad \quad \quad
  b_{21}^{(3)} = -\frac{
  b_{32}^{(3)} B_{33}^{(3)} +(B_{33}^{(3)})^2}{b_{32}^{(3)}}, \\
  b_{31}^{(3)} &= -b_{32}^{(3)} -B_{33}^{(3)}, \quad \quad \quad
  b_{32}^{(3)} = -\frac{2 (B_{33}^{(3)})^4 +1}{2 B_{33}^{(3)}
  ((B_{33}^{(3)})^2 +1)}.
\end{alignat*}
%
%
%
%
\section{Efficient drift-implicit SRK schemes and stability analysis}
\label{Sec:Optimal-Scheme-Stability-Ana}
The aim of this section is to determine efficient drift-implicit SRK
schemes that are included in the previously presented classification
with respect to a minimal number of implicit stages and explicit
function evaluations needed each step as well as good stability
properties. First, we briefly summarize the concept of mean--square
stability for SDEs. Therefore, we consider the scalar linear test
equation with multiplicative noise
\begin{equation} \label{Lin-stoch-test-eqn}
    {\mathrm{d}} X_t = \lambda \, X_t \, {\mathrm{d}}t + \mu \,
    X_t \, {\mathrm{d}}W_t ,
\end{equation}
for $t \geq t_0$ with initial value $X_{t_0}=x_0 \in \mathbb{R}
\setminus \{0\}$ and with some constants $\lambda,\mu \in
\mathbb{C}$. In order to analyse the mean--square stability
(MS--stability), we have to consider the second moment of the
solution process of SDE \eqref{Lin-stoch-test-eqn} and of the
corresponding numerical approximation process, respectively. The
solution of SDE \eqref{Lin-stoch-test-eqn} is said to be
(asymptotically) MS--stable if
\begin{equation} \label{Stability-domain-MS-mean-Loes}
    \lim_{t \to \infty} \E ( |X_t|^2 ) = 0 \quad
    \Leftrightarrow \quad 2 \, \Re(\lambda) + |\mu|^2 < 0
\end{equation}
holds for the coefficients $\lambda, \mu \in \mathbb{C}$, see e.~g.\
\cite{BuSi12,BuTi04,DeRoe08e,Hig00,Kloe,Sai} for further details. We call
$\mathcal{D}_{SDE} = \{ (\lambda,\mu) \in \mathbb{C}^2 : 2
\Re(\lambda)+|\mu|^2<0 \} \subset \mathbb{C}^2$ the domain of
MS--stability of SDE \eqref{Lin-stoch-test-eqn}. Here, we point out
that for $\mu=0$ the stability condition
\eqref{Stability-domain-MS-mean-Loes} reduces to the well known
deterministic stability condition $\Re(\lambda)<0$. \\ \\
In order to analyse the stability of the SRK method~\eqref{SSRK}, we
apply the method to the test problem \eqref{Lin-stoch-test-eqn}. We
are looking for conditions such that the SRK method yields
numerically stable solutions whenever
\eqref{Stability-domain-MS-mean-Loes} is fulfilled. A numerical
method is said to be numerically MS--stable if the approximations
$y_n$ satisfy $\lim_{n \to \infty} \E \left( |y_{n}|^2 \right) = 0$.
Applying the numerical method to \eqref{Lin-stoch-test-eqn}, we
obtain the recursion
\begin{equation} \label{Lemma-recursion-one-step-method}
    y_{n+1} = R_n(\hat{h}, k) \, y_n \, ,
\end{equation}
with a stability function $R_n(\hat{h}, k)$ using the
parametrization $\hat{h} = \lambda \, h$ and $k = \mu \sqrt{h}$ for
$h>0$ \cite{DeRoe08e,Hig00}.
Then, calculating the mean--square norm of
\eqref{Lemma-recursion-one-step-method}, we obviously yield
MS--stability, if
\begin{equation}
    \hat{R}(\hat{h},k):=\E(|R_n(\hat{h}, k)|^2)<1 .
\end{equation}
Now, we call $\mathcal{D}_{SRK} = \{ (\hat{h},k) \in \mathbb{C}^2 :
\hat{R}(\hat{h},k)<1 \} \subset \mathbb{C}^2$ the domain of
MS--stability of the SRK method. The numerical method is said to be
$A$--stable if $\mathcal{D}_{SDE} \subseteq \mathcal{D}_{SRK}$.
Because the domain of stability for $\lambda, \mu \in \mathbb{C}$ is
not easy to visualize, we have to restrict the figures to presenting
the region of stability for $\lambda, \mu \in \mathbb{R}$ in the
$\hat{h}$--$k^2$--plane. Then, for fixed values of $\lambda$ and
$\mu$, the set $\{(\lambda \, h, \mu^2 \, h) \subset \mathbb{R}^2 :
h > 0\}$ is a straight ray starting at the origin and going through
the point $(\lambda, \mu^2)$. Varying the step size $h$ corresponds
to moving along this ray. For $\lambda,\mu \in \mathbb{R}$, the
region of MS--stability for SDE \eqref{Lin-stoch-test-eqn} reduces
to the area of the $\hat{h}$--$k^2$--plane with the $\hat{h}$--axis
as the lower bound
and $k^2<-2 \hat{h}$ giving the upper bound for $\hat{h} < 0$. \\ \\
Next, we calculate the stability function $R_n(\hat{h},k)$ for the
$s$-stages SRK method \eqref{SSRK}. Let $H=(H_1,\dotsc,H_s)^T$.
Then \eqref{SSRK}
applied to \eqref{Lin-stoch-test-eqn} with equidistant step size
$h=h_n$ becomes
\[
    H = 
    e \, y_n + \lambda h A H + \mu \big(I_{(1),n}B^{(1)}
    + \frac{I_{(1,1),n}}{\sqrt{h}}B^{(2)} + \sqrt{h}B^{(3)}\big)H .
\]
Together with $I_{(1),n}=\sqrt{h_n} \xi_n$ where $\xi_n \sim N(0,1)$
and the parametrization $\hat{h}=\lambda h$ and $k=\mu\sqrt{h}$ this
can be reformulated to
\[
    H = \Big(I_s - \hat{h}A - k\big(\xi_n B^{(1)} + \frac{1}{2}(\xi_n^2-1)B^{(2)}
    + B^{(3)}\big) \Big)^{-1} e \, 
    y_n .
\]
Since the methods are stiffly accurate, that is $y_{n+1}=H_s$, the
stability function is given as
\begin{equation} \label{stabfun}
    R_n(\hat{h},k) = \varepsilon_s^T \Big(I_s - \hat{h}A - k\big(\xi_n
    B^{(1)} + \frac{1}{2}(\xi_n^2-1)B^{(2)} + B^{(3)}\big) \Big)^{-1} e 
\end{equation}
where $\varepsilon_s^T=(0,\dotsc,0,1)\in\mathbb{R}^s$.
\subsection{$A$-stable strong order 0.5 SRK schemes}
In the following, the computational costs are measured as the number
of function evaluations that are necessary in each step and we try
to minimize them. Therefore, the following coefficients for
drift-implicit order 0.5 SRK schemes are considered for both classes
I and II:
\begin{equation} \label{Sec:Opt-schemes-Order-0.5}
\begin{array}{|ccc|ccc|ccc|ccc|}
  \ \ a_1 \ \ && &&& &&&& 0 &&\\
  a_2 & 1-a_2 & \hspace*{5mm} & 1 & \hspace*{5mm} & \hspace*{5mm} & 0 &
  \hspace*{5mm} & \hspace*{5mm} & 0 & 0 & \hspace*{5mm} \\
  \hline
\end{array}
\end{equation}
where we choose $B_{11}^{(3)}=B_{21}^{(3)}=0$ and $a_1, a_2 \in
\mathbb{R}$.
\begin{figure}
\begin{center}
\includegraphics[width=6.8cm]{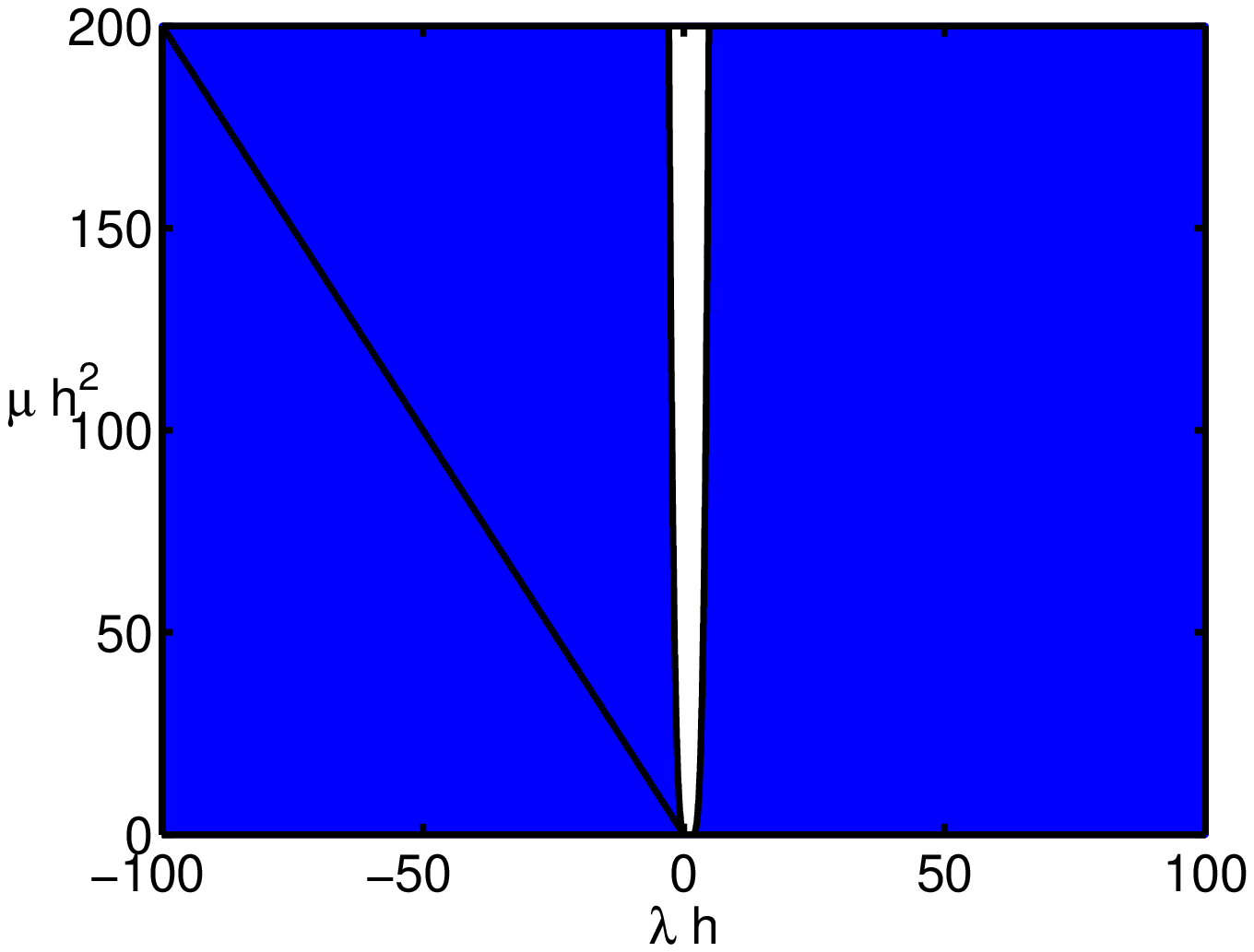}
\includegraphics[width=6.8cm]{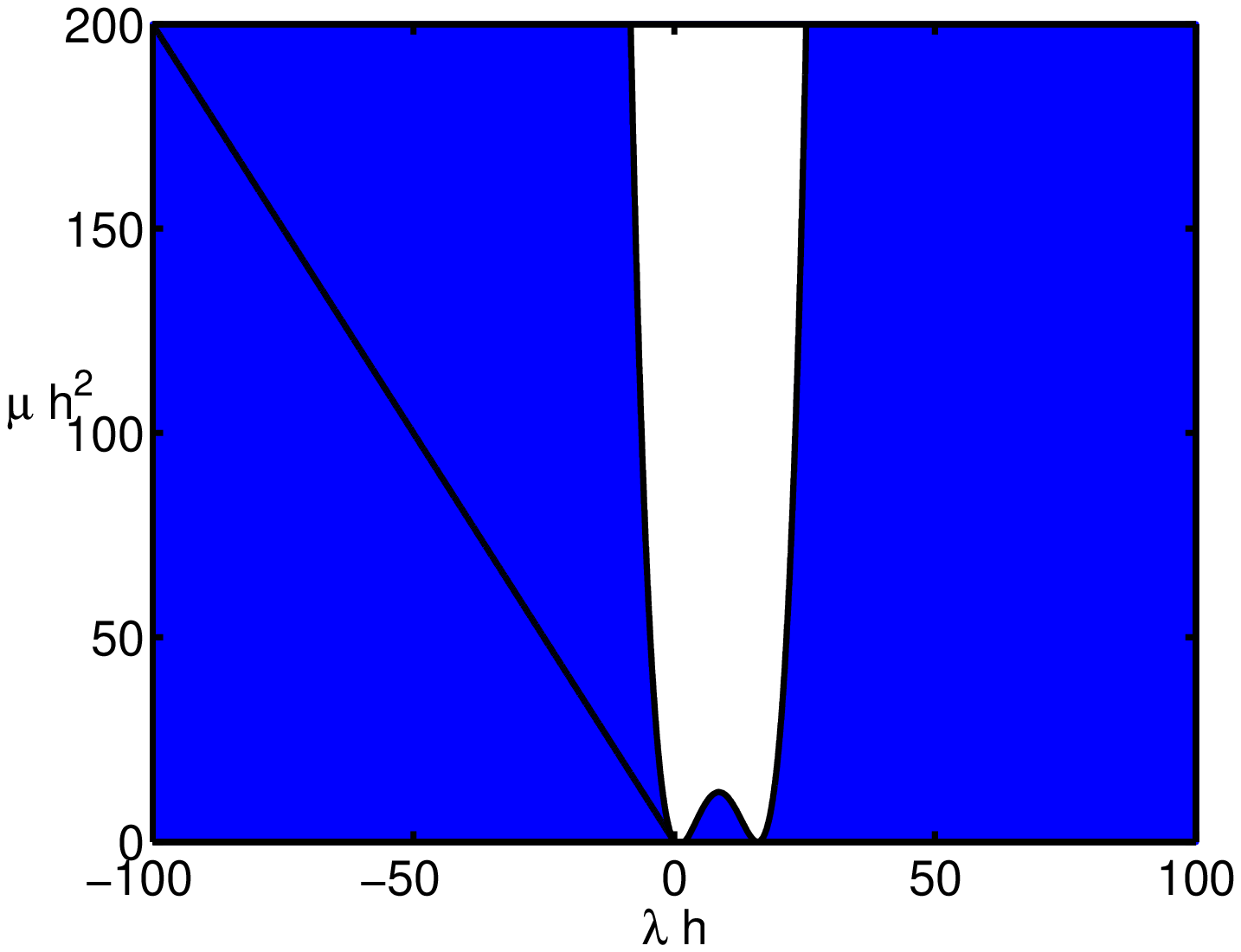}
\includegraphics[width=6.8cm]{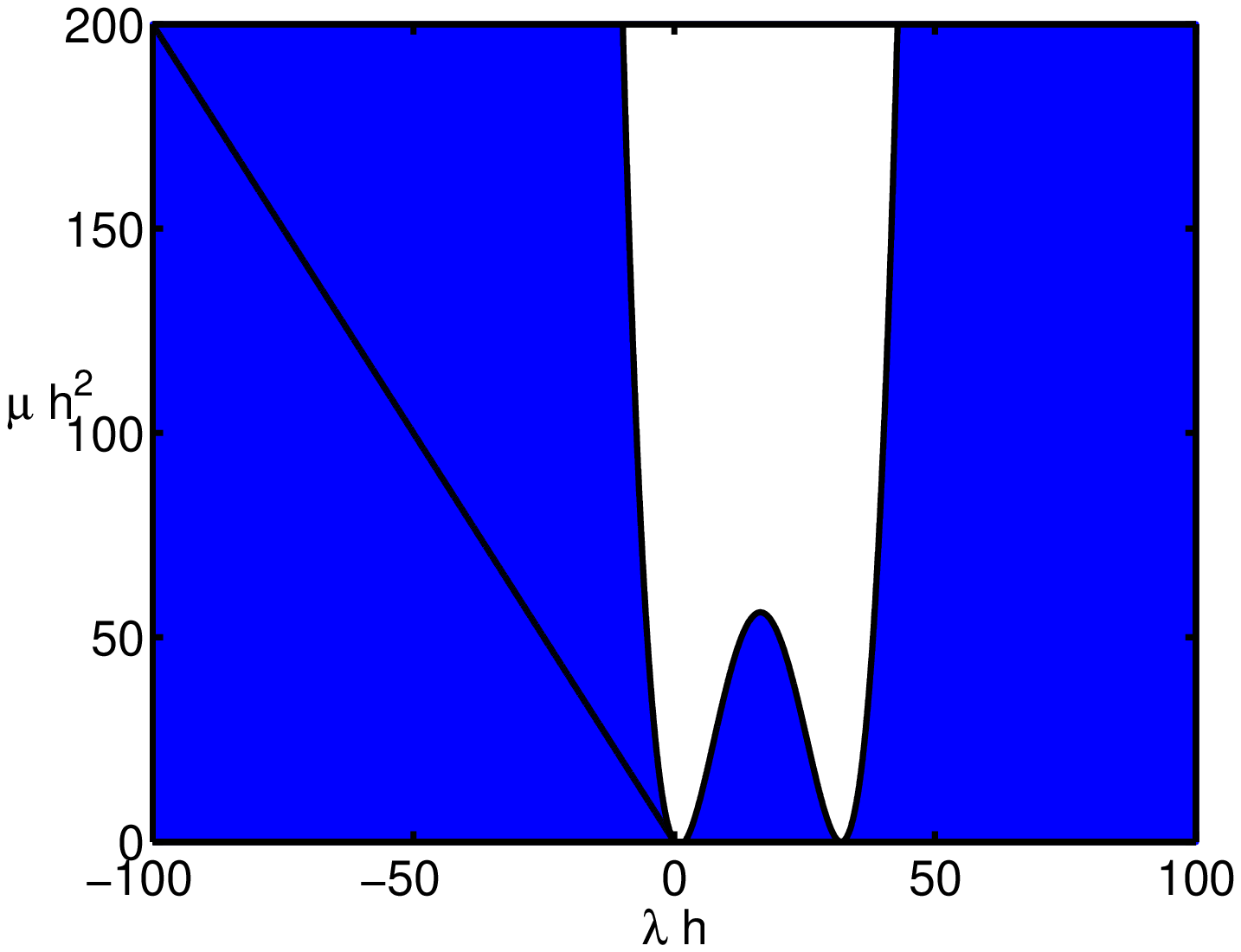}
\includegraphics[width=6.8cm]{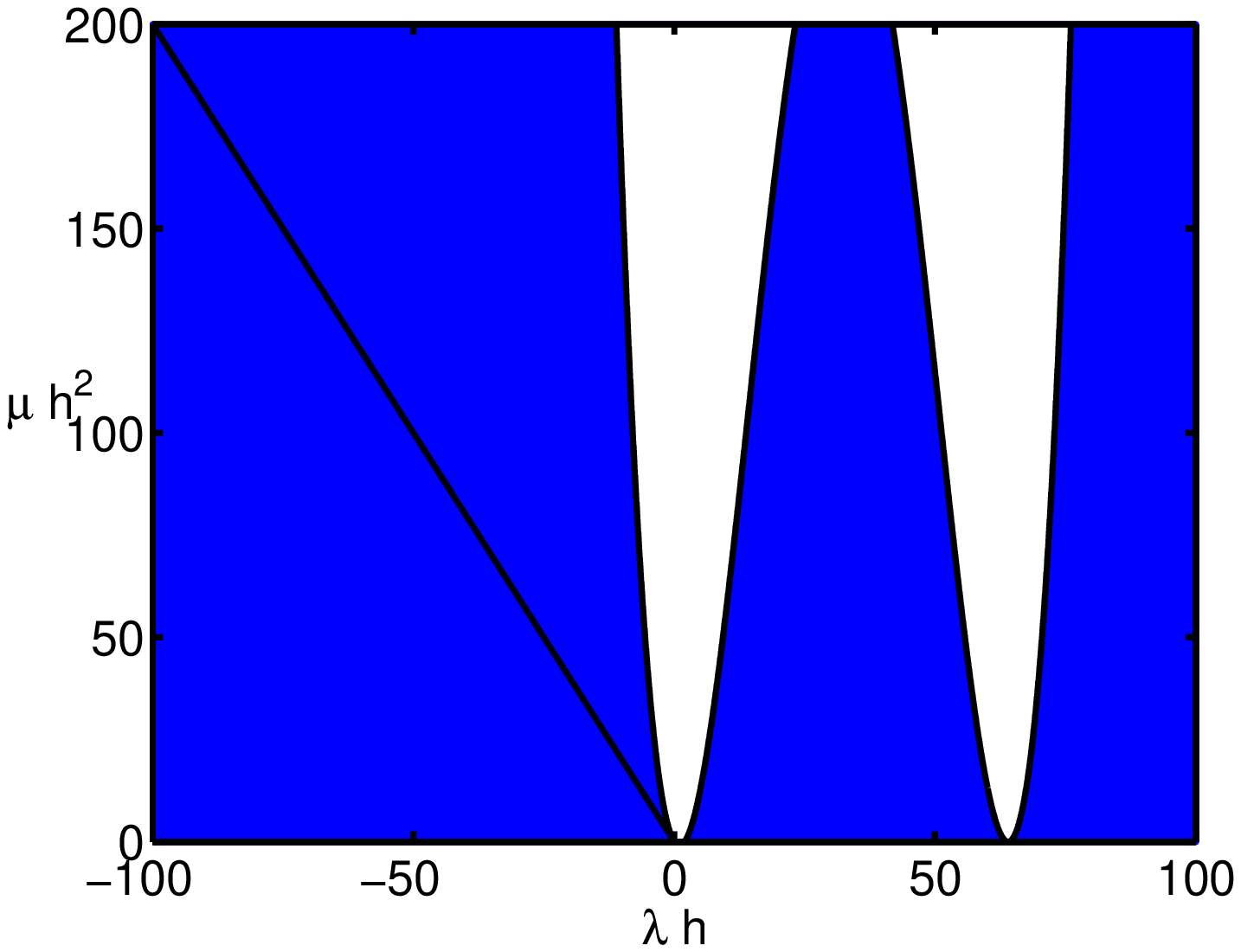}
\caption{Mean-square stability region for class I/II with $a_2=0$
and with $a_1=1$, $a_1=\frac{1}{16}$, $a_1=\frac{1}{32}$ and $a_1=\frac{1}{64}$, respectively.}
\label{MS-Bild-class-I-II-1a}
\end{center}
\end{figure}
First, we consider the case of diagonally drift-implicit SRK methods
where we choose $a_2=0$.
\begin{Lem} \label{Lemma-A-stab-Order-0.5}
    The order 0.5 SRK scheme with coefficients
    \eqref{Sec:Opt-schemes-Order-0.5} and $a_2=0$ is $A$-stable for equation
    \eqref{Lin-stoch-test-eqn}, i.e.\ $\mathcal{D}_{SDE} \subseteq
    \mathcal{D}_{SRK}$, if and only if $a_1 \geq 0$.
\end{Lem}
{\textbf{Proof.}} Calculating $\hat{R}(\hat{h},k)$ from the
stability function \eqref{stabfun}
using the coefficients \eqref{Sec:Opt-schemes-Order-0.5} yields
\begin{equation} \label{Proof-0.5-eqn-1}
    \hat{R}(\hat{h},k) = \frac{|1-a_1 \hat{h}|^2 + |k|^2}{|1-\hat{h}|^2
    |1-a_1 \hat{h}|^2} \, .
\end{equation}
Now, we obtain that $\mathcal{D}_{SDE} \subseteq \mathcal{D}_{SRK}$
if $\hat{R}(\hat{h},k) < 1$ for all $\hat{h}, k \in \mathbb{C}^2$
with $2 \Re(\hat{h}) + |k|^2 <0$. Assuming that $|k|^2 < -2
\Re(\hat{h})$, we have to prove that $\hat{R}(\hat{h},k)-1 < 0$
holds. Using this assumption, we get
\begin{equation} \label{Proof-0.5-eqn-2}
    \hat{R}(\hat{h},k)-1 < \frac{\phi(\hat{h},a_1)}{|1-\hat{h}|^2
    |1-a_1 \hat{h}|^2}
\end{equation}
with
\begin{equation} \label{Proof-0.5-eqn-3}
    \begin{split}
        \phi(\hat{h},a_1) :=  &\,(-4a_1 -1) \Re(\hat{h})^2 + (2a_1^2+2a_1)
        \Re(\hat{h})^3 -a_1^2  \Re(\hat{h})^4 -\Im(\hat{h})^2 \\
        &+(2a_1^2 + 2a_1)  \Re(\hat{h})  \Im(\hat{h})^2 - 2a_1^2  \Re(\hat{h})^2
        \Im(\hat{h})^2 - a_1^2  \Im(\hat{h})^4 \, .
    \end{split}
\end{equation}
Now, for $\Re(\hat{h})<0$ and $a_1 \geq 0$ the expression
\eqref{Proof-0.5-eqn-3} is obviously not positive, i.e., the order
0.5 scheme \eqref{Sec:Opt-schemes-Order-0.5} is $A$-stable. \\ \\
For $a_1<0$, we restrict our analysis to the case of
$\Im(\hat{h})=\Im(k)=0$ in the following. Since $\hat{R}(\hat{h},k)$
has a singularity at $\hat{h}=\tfrac{1}{a_1}$, we restrict our
considerations to the case where $\hat{h}<\tfrac{1}{a_1}$. Then,
considering the boundary $|k|=\sqrt{-2 \hat{h}}$ of the domain of
stability of the test equation \eqref{Lin-stoch-test-eqn}, we get
\begin{equation} \label{Proof-0.5-eqn-4}
    \hat{R}(\hat{h},\sqrt{-2 \hat{h}}) -1 = \frac{\hat{h}^2 \left(
    \hat{h}^2 -\frac{2+2a_1}{a_1} \hat{h} + \frac{4a_1+1}{a_1^2}
    \right)}{|1-\hat{h}|^2 |1-a_1 \hat{h}|^2} \, .
\end{equation}
By calculating the roots of \eqref{Proof-0.5-eqn-4} we get that
$\hat{R}(\hat{h}, \sqrt{-2 \hat{h}})-1>0$ for
\begin{equation} \label{Proof-0.5-eqn-5}
    \hat{h} \in \, I(a_1) := \left] \tfrac{1+a_1+\sqrt{a_1^2-2a_1}}{a_1},
    \tfrac{1}{a_1} \right[ \, .
\end{equation}
Due to the continuity of $\hat{R}(\hat{h}, k)$ on
$\,]-\infty,\tfrac{1}{a_1}[\, \times \mathbb{R}$ there exists some
$\varepsilon>0$ such that $\hat{R}(\hat{h},k)-1>0$ on some open ball
$B_{\varepsilon}(\hat{h},\sqrt{-2 \hat{h}})$ with radius
$\varepsilon$ and center $(\hat{h},\sqrt{-2 \hat{h}})$ with $\hat{h}
\in I(a_1)$. Since $B_{\varepsilon}(\hat{h},\sqrt{-2 \hat{h}}) \cap
\mathcal{D}_{SDE} \neq \emptyset$ it follows that the scheme can not
be $A$-stable. \hfill $\Box$ \\ \\
Considering the regions of MS-stability for the SRK schemes with
$a_2=0$ and different values $a_1 \in \{1, \tfrac{1}{16},
\tfrac{1}{32}, \tfrac{1}{64}\}$, we can see in
Figure~\ref{MS-Bild-class-I-II-1a} that $\mathcal{D}_{SDE} \subseteq
\mathcal{D}_{SRK}$ is always fulfilled.
\begin{Bem}
    If we choose $a_1=1$ and $a_2=0$ in
    \eqref{Sec:Opt-schemes-Order-0.5}, then the resulting order 0.5
    scheme is $A$-stable and a singly diagonally drift-implicit
    stiffly accurate SRK scheme.
    Especially, the calculation of only one $LU$
    decomposition is needed each step if a simplified Newton method is applied
    to solve the implicit equations.
\end{Bem}
As another class of schemes, we consider the case of an explicit
first stage, i.e.\ where $a_1=0$. However, then we need $a_2 \neq 1$
if the SRK method is applied to an SDAE, see \cite{KKR12}.
\begin{figure}
\begin{center}
\includegraphics[width=6.8cm]{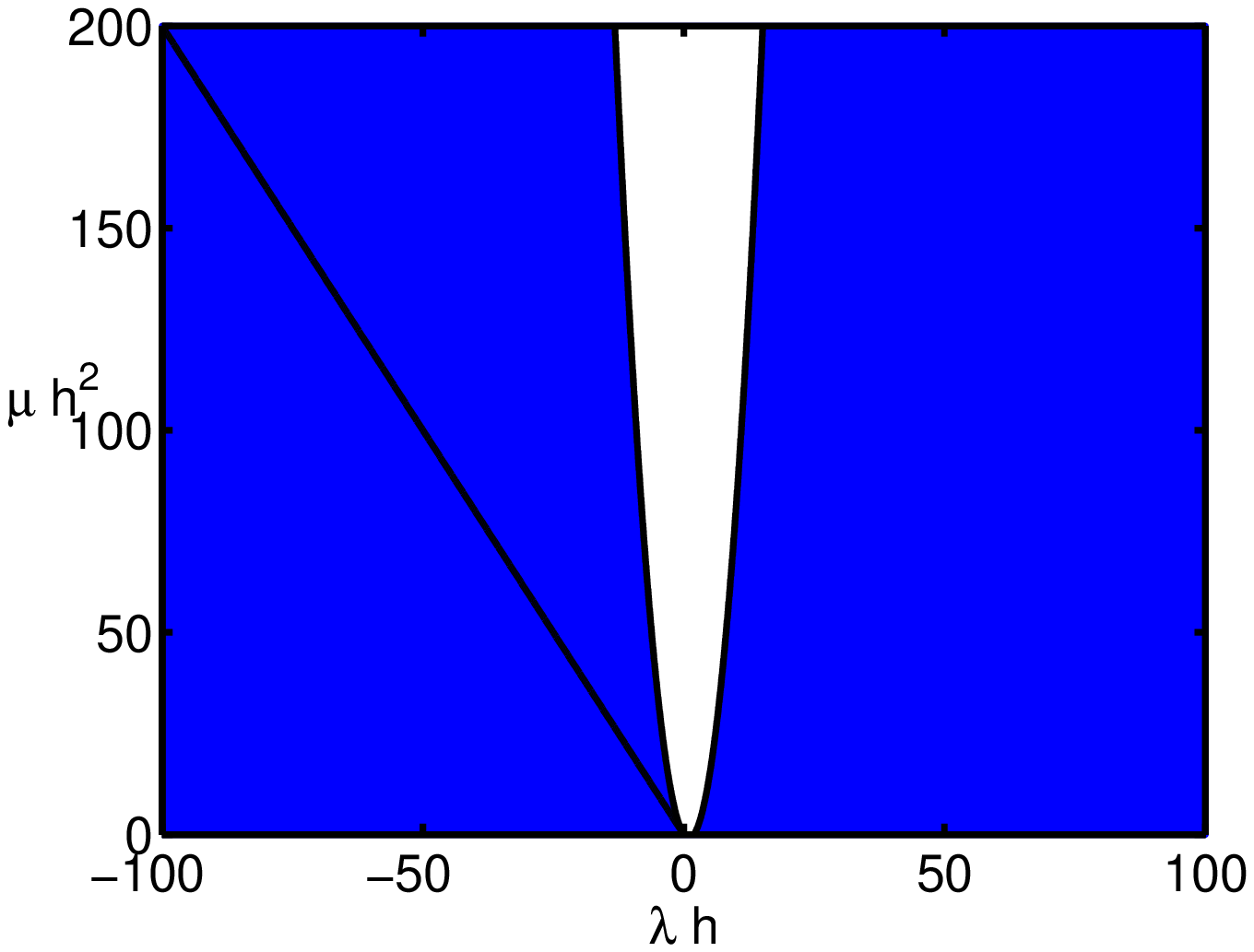}
\includegraphics[width=6.8cm]{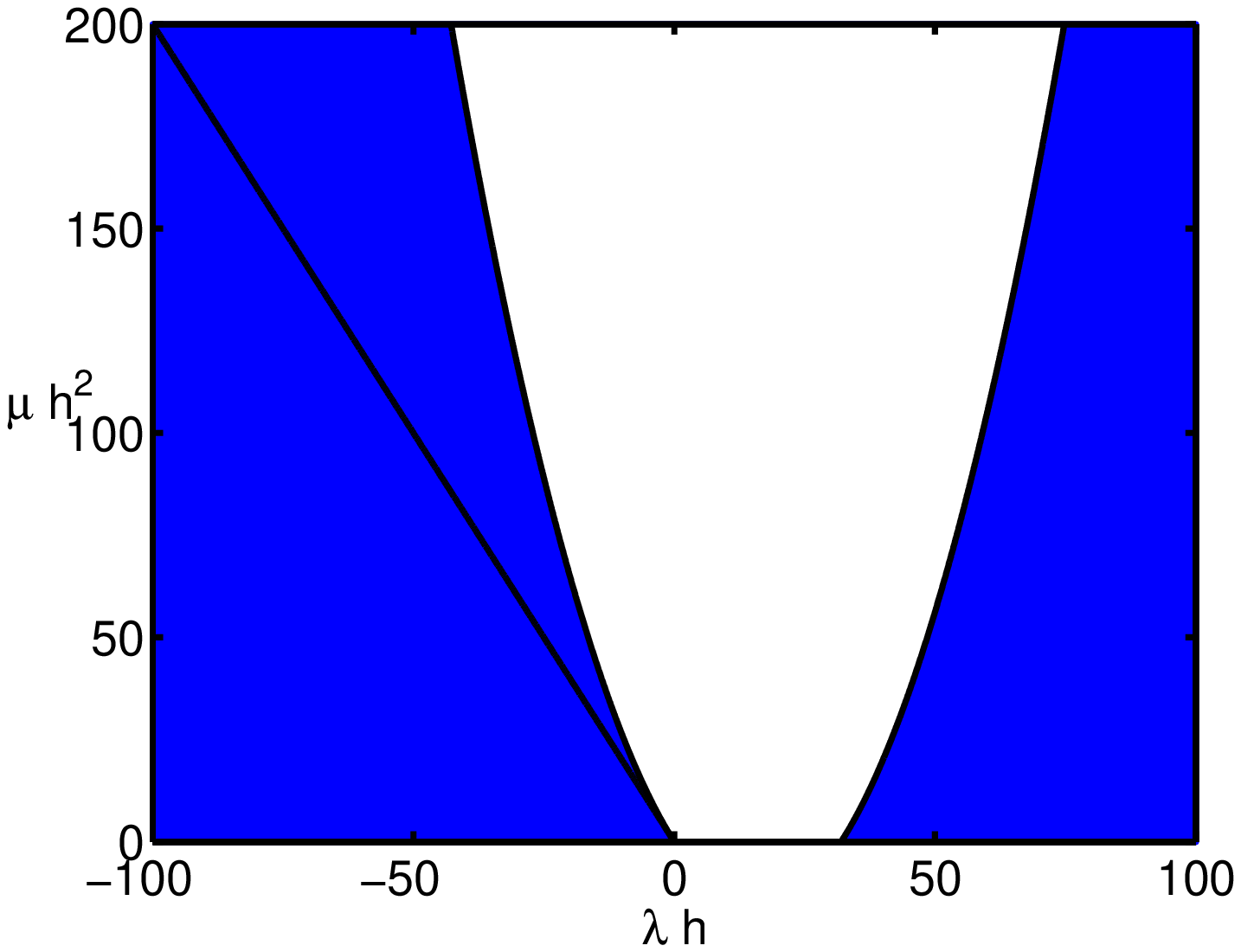}
\includegraphics[width=6.8cm]{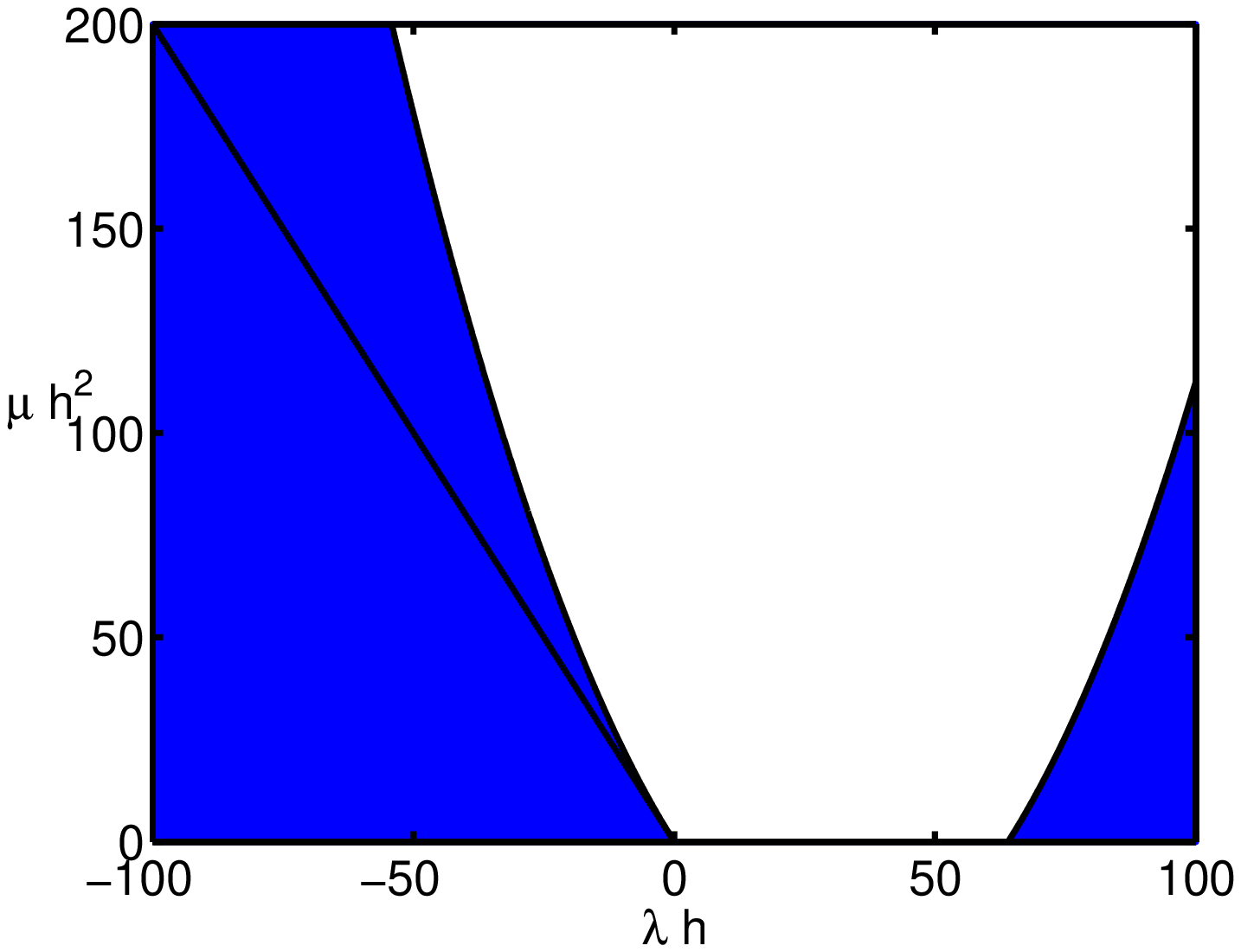}
\includegraphics[width=6.8cm]{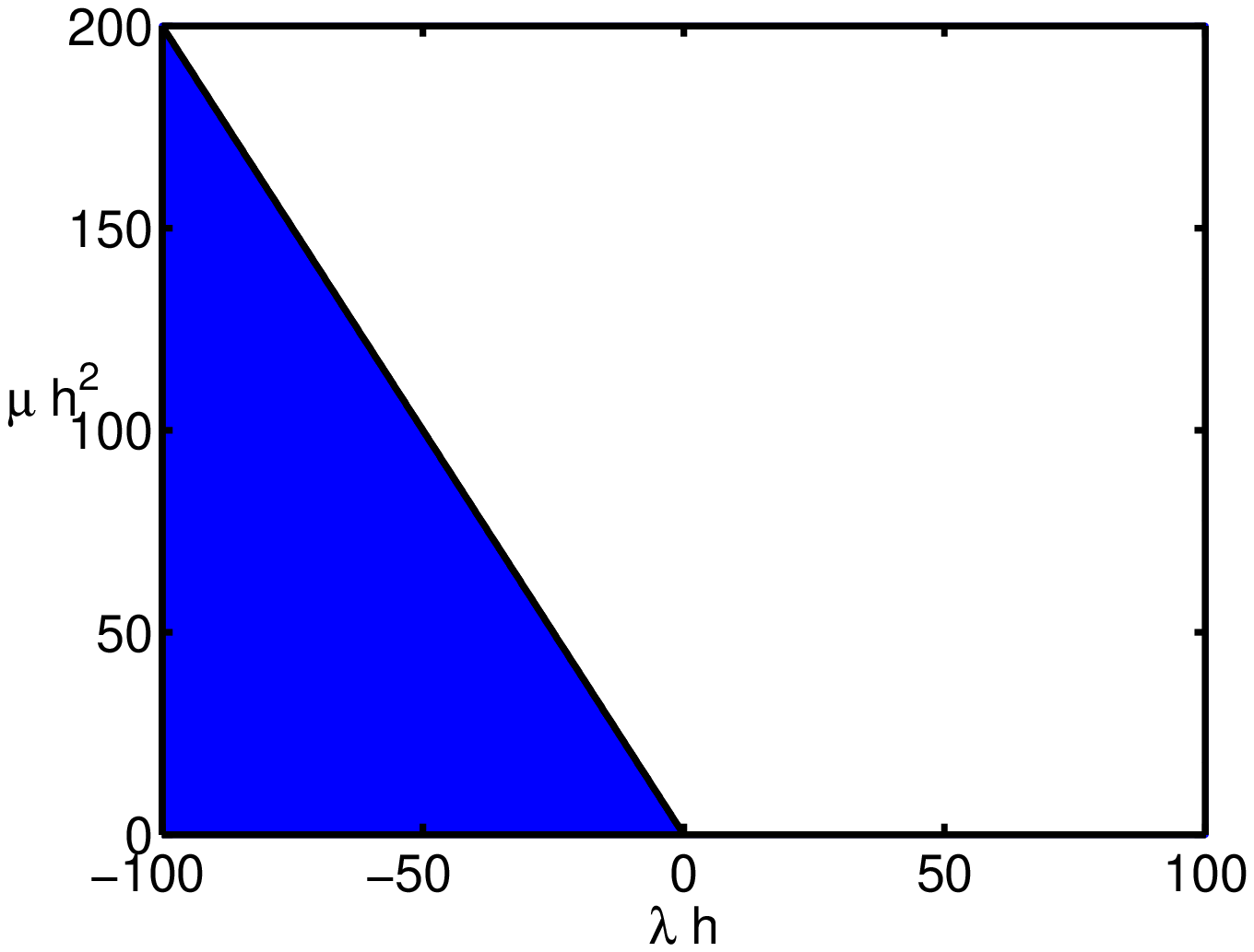}
\caption{Mean-square stability region for class I/II with $a_1=0$
and $a_2=0$, $a_2=\frac{15}{32}$, $a_2=\frac{31}{64}$ and with $a_2=\tfrac{1}{2}$,
respectively.} \label{MS-Bild-class-I-II-1b}
\end{center}
\end{figure}
\begin{Lem} \label{Lemma-A-stab-Order-0.5-Allg}
    The order 0.5 SRK scheme with coefficients
    \eqref{Sec:Opt-schemes-Order-0.5} is $A$-stable for equation
    \eqref{Lin-stoch-test-eqn}, i.e.\ $\mathcal{D}_{SDE} \subseteq
    \mathcal{D}_{SRK}$, if $a_1 \geq 0$ and $a_2 \leq \min \{
    \frac{1+4 a_1}{2(1+a_1)}, \frac{1}{2(1-a_1)}, 1 \}$.
\end{Lem}
{\textbf{Proof.}} Calculating $\hat{R}(\hat{h},k)$ from the
stability function \eqref{stabfun}
using the coefficients \eqref{Sec:Opt-schemes-Order-0.5} yields
\begin{equation} \label{Proof-0.5-eqn-1-A}
    \hat{R}(\hat{h},k) = \frac{|1+(a_2-a_1) \hat{h}|^2 + |k|^2}{|1- (1-a_2) \hat{h}|^2
    |1-a_1 \hat{h}|^2} \, .
\end{equation}
Now, we obtain that $\mathcal{D}_{SDE} \subseteq \mathcal{D}_{SRK}$
if $\hat{R}(\hat{h},k) < 1$ for all $\hat{h}, k \in \mathbb{C}^2$
with $2 \Re(\hat{h}) + |k|^2 <0$. Assuming that $|k|^2 < -2
\Re(\hat{h})$, we have to prove that $\hat{R}(\hat{h},k)-1 < 0$
holds. Using this assumption, we get
\begin{equation} \label{Proof-0.5-eqn-2-A}
    \hat{R}(\hat{h},k)-1 < \frac{\phi(\hat{h},a_1,a_2)}{|1-(1-a_2) \hat{h}|^2
    |1-a_1 \hat{h}|^2}
\end{equation}
with
\begin{equation} \label{Proof-0.5-eqn-3-A}
    \begin{split}
        \phi(\hat{h},a_1,a_2) :=  &\,(-4a_1 -1+2a_2+2a_1 a_2) \Re(\hat{h})^2
        + (2a_2-2a_1 a_2-1) \Im(\hat{h})^2 \\
        &+ 2(a_1^2(1-a_2)+a_1 (1-a_2)^2) |\hat{h}|^2 \Re(\hat{h})
        -a_1^2 (1-a_2)^2 |\hat{h}|^4 \, .
    \end{split}
\end{equation}
Thus, for $\Re(\hat{h})<0$ the expression \eqref{Proof-0.5-eqn-3-A}
is obviously not positive if $a_1 \geq 0$ and if
\begin{equation*}
    a_2 \leq \min \big\{ \frac{1+4 a_1}{2(1+a_1)}, \frac{1}{2(1-a_1)}, 1 \big\}
    \, .
\end{equation*}
Then, the order 0.5 scheme \eqref{Sec:Opt-schemes-Order-0.5} is
$A$-stable. \hfill $\Box$ \\ \\
In case of $a_1=0$, the regions of MS-stability for the SRK method
with $a_2 \in \{0, \tfrac{15}{32}, \tfrac{31}{64}, \tfrac{1}{2}\}$
are presented in Figure~\ref{MS-Bild-class-I-II-1b} where
$\mathcal{D}_{SDE} \subseteq \mathcal{D}_{SRK}$ is fulfilled.
For $a_1=0$ and $a_2=\frac{1}{2}$, the region of MS-stability for
the SRK scheme coincides perfectly with the region of MS-stability
for the test SDE.
\begin{Bem}
    In the case of $a_1=0$ and $a_2 \neq 1$ in \eqref{Sec:Opt-schemes-Order-0.5}, the
    order 0.5 stiffly accurate drift-implicit SRK scheme coincides
    with the well known $\theta$-method \cite{Hig00} and
    needs only one
    stage-evaluation of the drift function $f$ and one of the diffusion
    function $g$ each step due to the FSAL (first same as last) property
    \cite{Hai2}. Further, only one implicit equation has to be
    solved each step.
\end{Bem}
\subsection{$A$-stable strong order 1.0 SRK schemes}
Next, we want to find some $A$-stable order 1.0 SRK schemes.
As mentioned in Section~\ref{Sec:Class-order-1.0} the smallest
number of stages for order 1.0 schemes is $s=3$. Within this case of
$3$-stages schemes, it turns out that the Classes II and X are the
ones with the lowest number of function evaluations, i.e.\ with
minimal computational costs. This is due to the fact that these are
the classes including schemes that are explicit in the diffusion.
\\ \\
Thus, choosing the coefficients for Class II such that the
computational effort is minimized, i.e.\ with $A_{11}=a_1$,
$A_{22}=a_2$, $A_{33}=a_3$, $B_{22}^{(3)}=0$ and $B_{32}^{(3)}=\pm
\frac{1}{2b}$, we get the tableau
\begin{equation} \label{Sec:Opt-schemes-Class-II}
    \begin{tabular}{|ccc|ccc|ccc|ccc|}
     $a_1$ & & &\hspace*{5mm}&\hspace*{5mm}& \hspace*{5mm}
     &\hspace*{5mm}&\hspace*{5mm}&\hspace*{5mm}& $0$ &&\\
     $a_1-a_2$ & $a_2$ & & $b$ & & & $0$ & && $\mp b$
     & $0$ & \\
     $1-a_3$ & $0$ & $a_3$ & $1-\tfrac{1}{2b}$ & $\tfrac{1}{2b}$& &
     $0$ & $0$ & & $\mp \tfrac{1}{2b}$ & $\pm \tfrac{1}{2 b}$ &
     $0$ \\ \hline
   \end{tabular}
\end{equation}
with $a_1, a_2, a_3 \in \mathbb{R}$ and $b \in \mathbb{R} \setminus
\{0\}$. \\ \\
Further, choosing the coefficients for Class X such that the
computational effort is minimized, i.e.\ with $A_{11}=a_1$,
$A_{22}=a_2$, $A_{33}=a_3$, $A_{21}=a_4$,
$B_{22}^{(3)}=0$ and $B_{32}^{(2)}=\frac{1}{b}$ results
in the tableau
\begin{equation} \label{Sec:Opt-schemes-Class-X}
    \begin{tabular}{|ccc|ccc|ccc|ccc|}
     $a_1$ & & &\hspace*{5mm}&\hspace*{5mm}& \hspace*{5mm}
     &\hspace*{5mm}&\hspace*{5mm}&\hspace*{5mm}& $0$ &&\\
     $a_4$ & $a_2$ & & $0$ & & & $0$ & && $b$
     & $0$ & \\
     $1-a_3$ & $0$ & $a_3$ & $1$ & $0$& &
     $-\frac{1}{b}$ & $\frac{1}{b}$ & & $0$ & $0$ &
     $0$ \\ \hline
   \end{tabular}
\end{equation}
with $a_1, a_2, a_3, a_4 \in \mathbb{R}$ and $b \in \mathbb{R} \setminus
\{0\}$. \\ \\
Now, in the case of $s=3$ with $A_{ij}=B^{(3)}_{ij}=0$ for $j>i$ and
$B^{(1)}_{ij}=B^{(2)}_{ij}=0$ for $j \geq i$, by rearranging the
terms with respect to powers of $\xi_n$ the stability function
\eqref{stabfun} has a representation of type
\begin{equation*} 
    R_n(\hat{h},k) = \Gamma + \Sigma_1 \xi_n + \Sigma_2 \xi_n^2 +
    \Sigma_3 \xi_n^3 + \Sigma_4 \xi_n^4
\end{equation*}
with some suitable coefficients $\Gamma, \Sigma_1, \ldots, \Sigma_4$
independent of $\xi_n$, see also \cite{KKR12}.
Therefore, we calculate the mean-square stability function
$\hat{R}(\hat{h},k)$ for the diagonally implicit SRK method
\eqref{SSRK} as
\begin{equation} \label{Stab-func}
  \begin{split}
  \hat{R}(\hat{h},k) = &|\Gamma|^2 + \Gamma \overline{\Sigma}_2 +
  \overline{\Gamma} \Sigma_2 + 3 \Gamma \overline{\Sigma}_4 + 3
  \overline{\Gamma} \Sigma_4 + |\Sigma_1|^2 + 3 \Sigma_1
  \overline{\Sigma}_3 + 3 \overline{\Sigma}_1 \Sigma_3 \\
  &+ 3 |\Sigma_2|^2 + 15 \Sigma_2 \overline{\Sigma}_4 + 15
  \overline{\Sigma}_2 \Sigma_4 + 15 |\Sigma_3|^2 + 105 |\Sigma_4|^2
  \, .
  \end{split}
\end{equation}
Especially, for class II with the coefficients
\eqref{Sec:Opt-schemes-Class-II}, we get
\begin{equation} \label{Class-II-Gamma}
    \Gamma = \frac{1-\tfrac{1}{2} k^2 - (a_1+a_2+a_3-1) \hat{h} + a_2(a_1+a_3-1)
    \hat{h}^2}{(1-a_1 \hat{h}) (1-a_2 \hat{h}) (1-a_3 \hat{h})} ,
\end{equation}
\begin{equation} \label{Class-II-Sigma-1-2}
    \Sigma_1 = \frac{k}{(1-a_1 \hat{h})(1-a_3 \hat{h})} , \qquad
    \Sigma_2 = \frac{k^2}{2 (1-a_1 \hat{h})(1-a_2 \hat{h}) (1-a_3 \hat{h})} ,
\end{equation}
and $\Sigma_3 = \Sigma_4 = 0$. Here, we would like to point out,
that the stability function does not depend on the parameter $b$.
Further, for class X with the coefficients
\eqref{Sec:Opt-schemes-Class-X}, we have
\begin{equation} \label{Class-X-Gamma}
    \Gamma = \frac{1-\tfrac{1}{2} k^2 + (1-a_1-a_2-a_3) \hat{h} + (a_1+a_3-1) a_2
    \hat{h}^2 + (a_1-a_2-a_4) \frac{1}{2b} k \hat{h}}{(1-a_1
    \hat{h}) (1-a_2 \hat{h}) (1-a_3 \hat{h})} ,
\end{equation}
\begin{equation} \label{Class-X-Sigma-1-2}
    \Sigma_1 = \frac{k}{(1-a_1 \hat{h})(1-a_3 \hat{h})} , \qquad
    \Sigma_2 = \frac{\frac{1}{2}k^2+(a_2-a_1+a_4) \frac{1}{2b} k
    \hat{h}}{(1-a_1 \hat{h}) (1-a_2 \hat{h}) (1-a_3 \hat{h})} ,
\end{equation}
and $\Sigma_3 = \Sigma_4 = 0$. \\ \\
First, we consider the case of diagonally drift-implicit stiffly
accurate SRK methods. Therefore, we analyse class II with
$a_1=a_2=a$ for some $a \in \mathbb{R}$ and $a_3=1$. Here, we have
to point out, that we need $a \neq 0$ if the SRK method is applied
to SDAEs, see \cite{KKR12}. Then, we need three stage-evaluations of
the drift function $f$ and two stage-evaluations of the diffusion
function $g$ for the diagonally implicit SRK method \eqref{SSRK}
each step.
\begin{figure}
\begin{center}
\includegraphics[width=6.8cm]{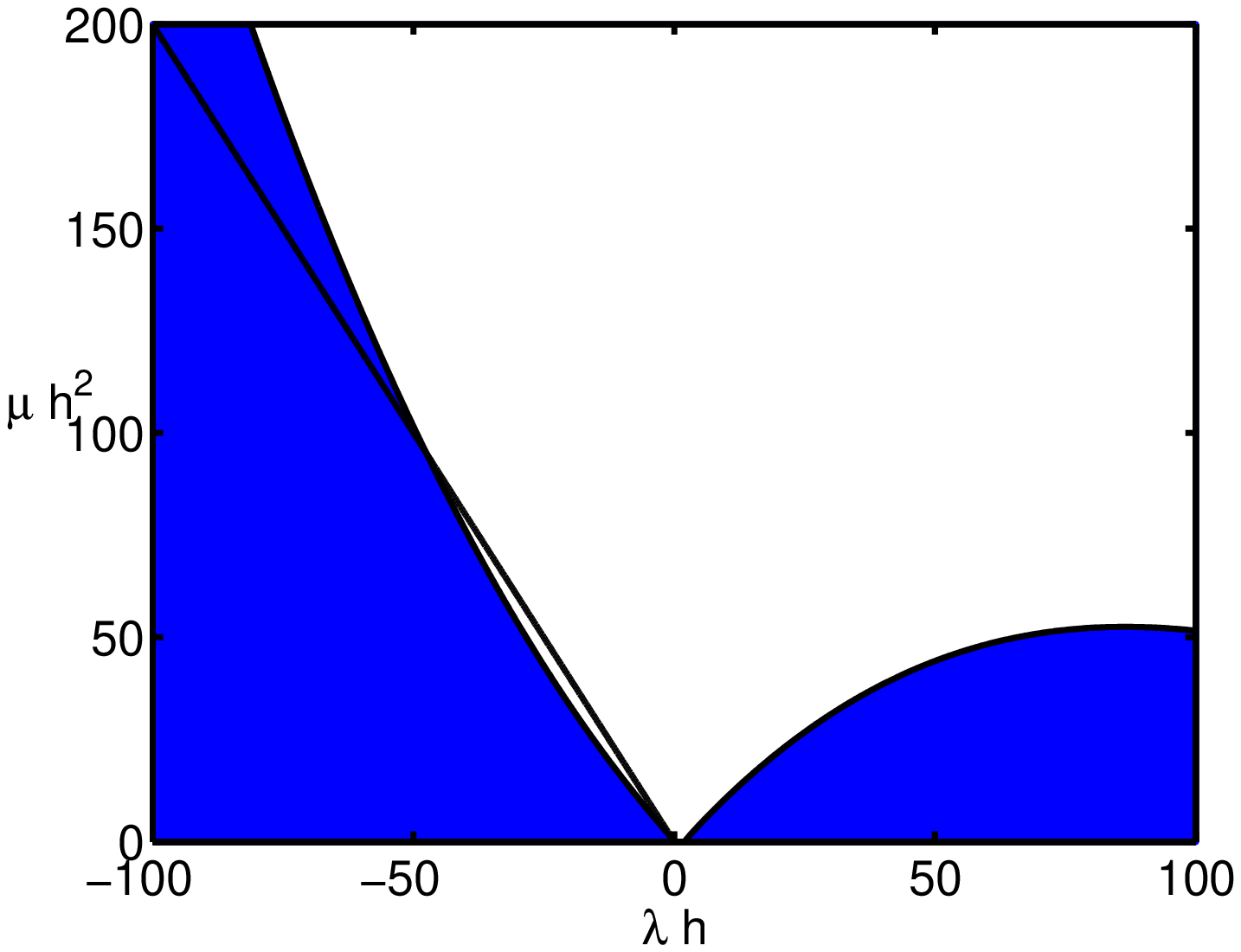}
\includegraphics[width=6.8cm]{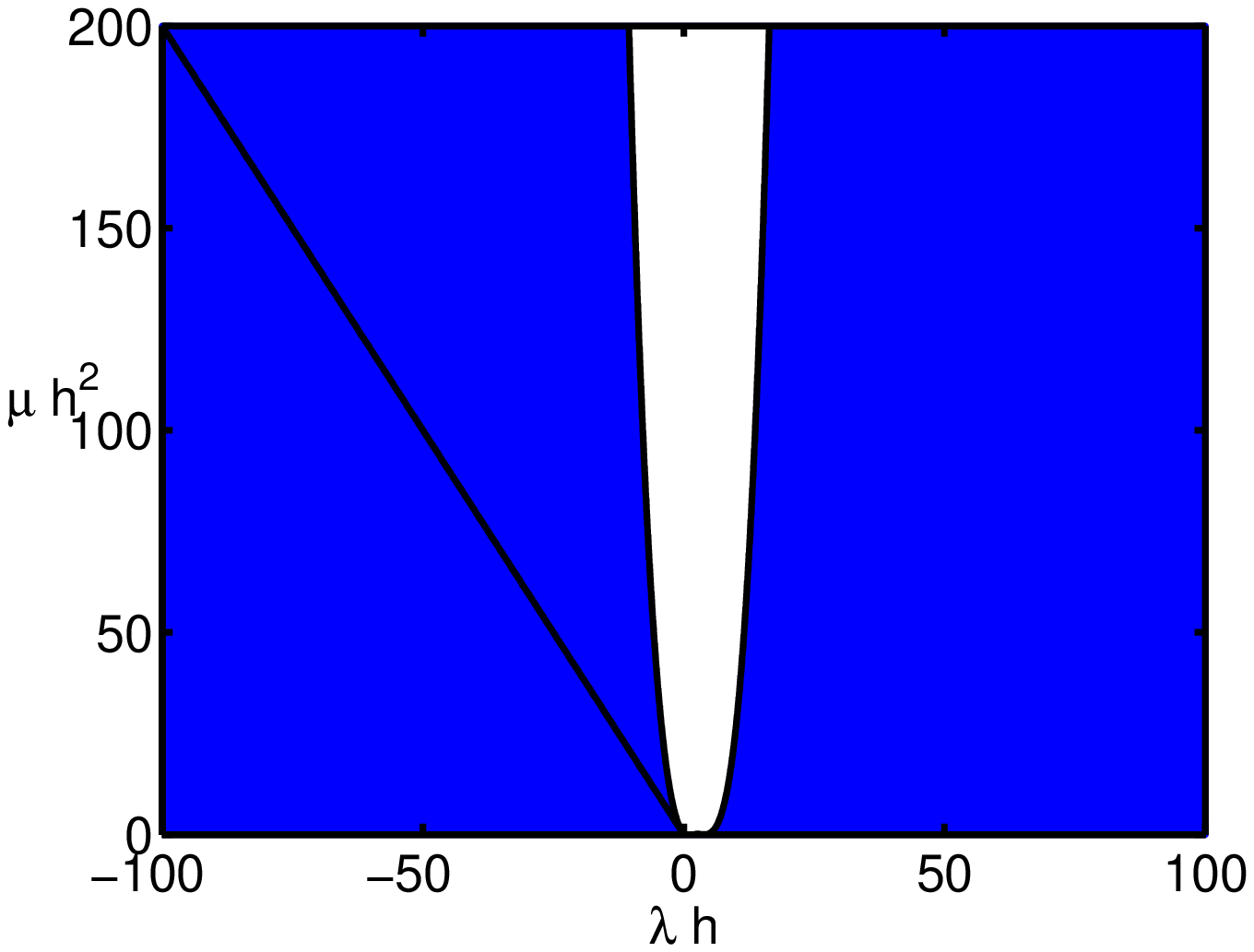}
\includegraphics[width=6.8cm]{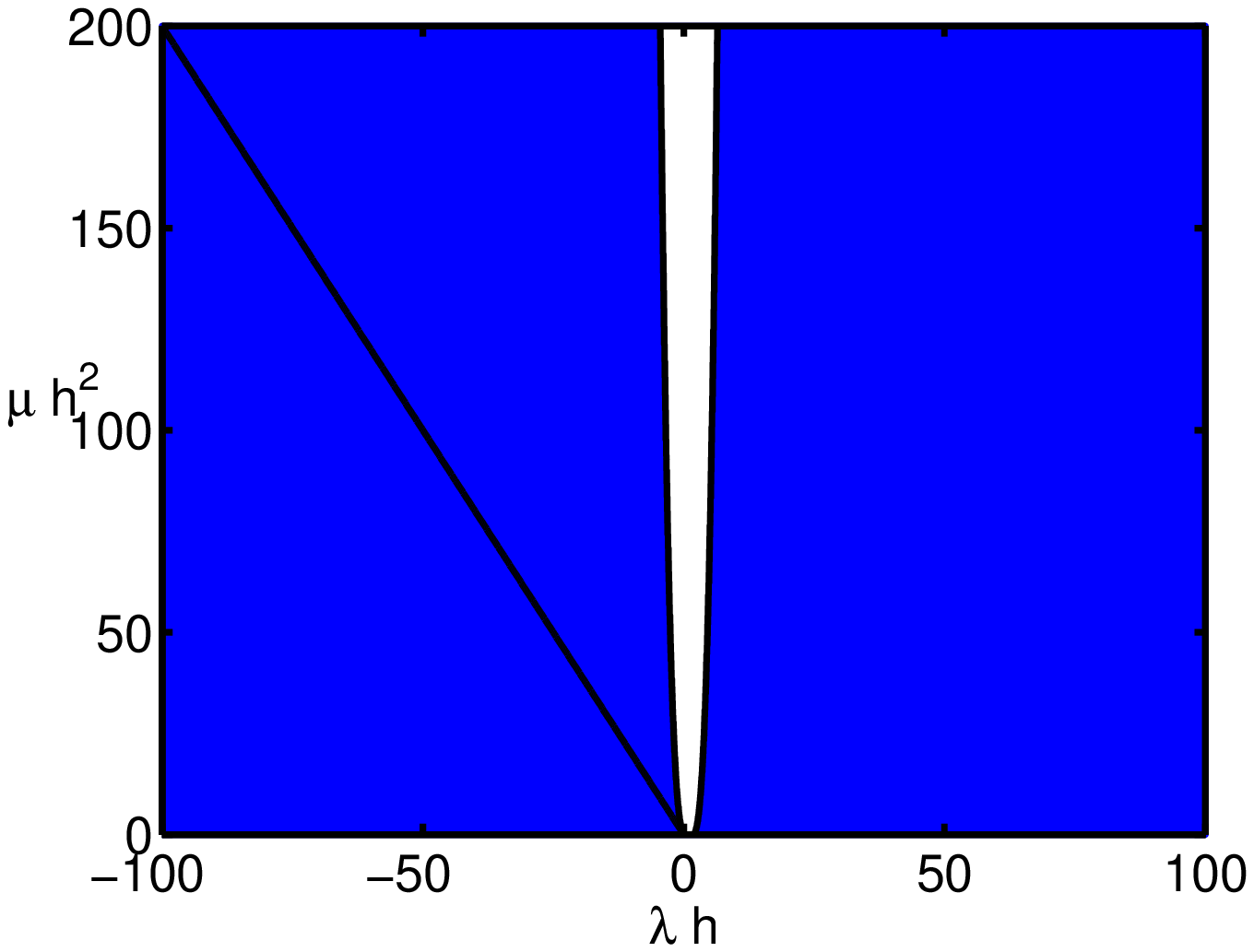}
\includegraphics[width=6.8cm]{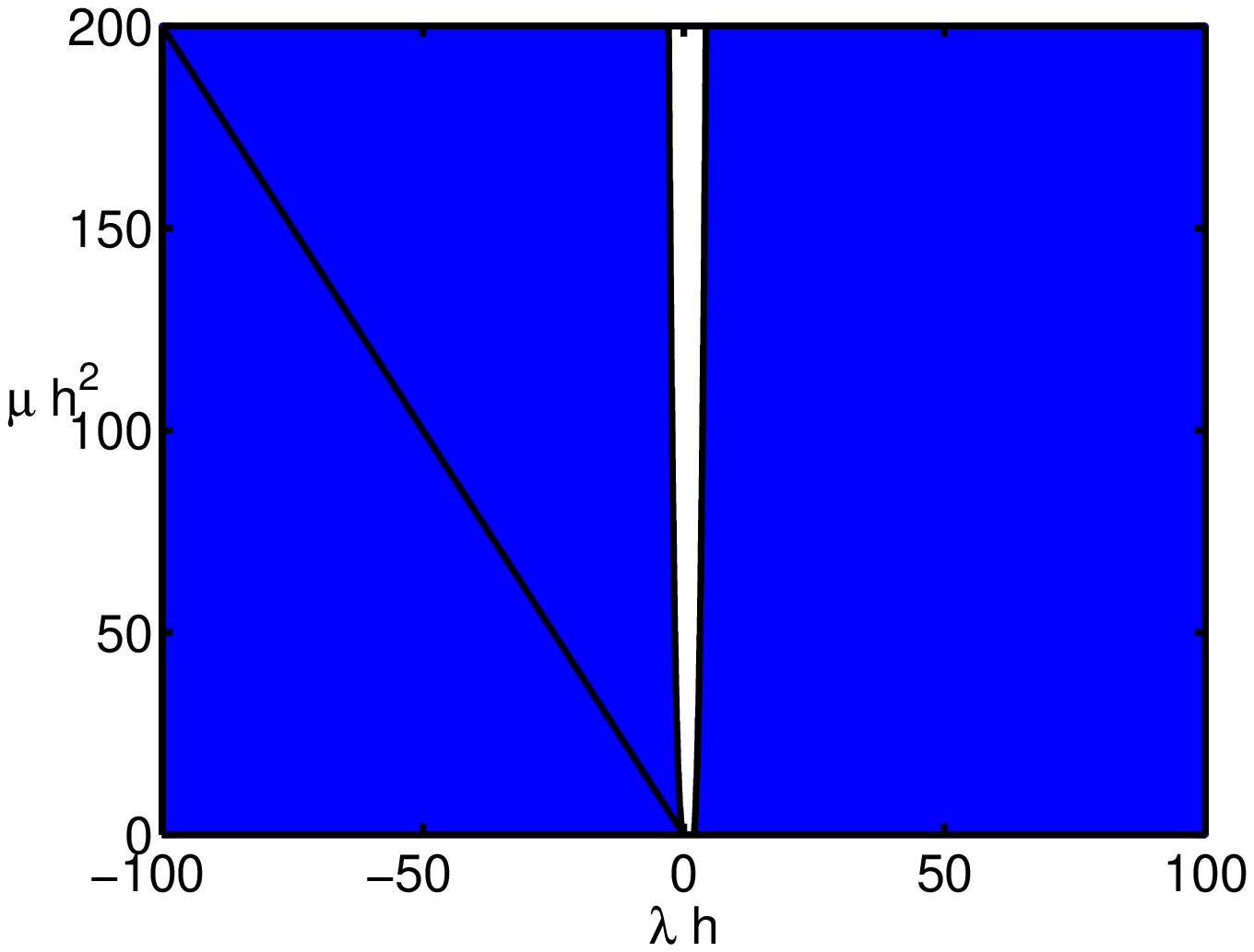}
\caption{Mean-square stability region for class II and class X with $a_3=1$, $a_4=0$, $b=1$ and with
$a_1=a_2=\frac{1}{256}$ (not $A$-stable), $a_1=a_2=\frac{1}{4}$, $a_1=a_2=1$ and $a_1=a_2=2$, respectively.}
\label{MS-Bild-class-II-X-1}
\end{center}
\end{figure}
\begin{Lem} \label{Lemma-A-stab-Class-II}
  The family of order 1.0 SRK schemes with coefficients \eqref{Sec:Opt-schemes-Class-II}
  in case of $a_1=a_2=a$ and $a_3=1$ is $A$-stable for equation \eqref{Lin-stoch-test-eqn},
  i.e.\ it holds $\mathcal{D}_{SDE} \subseteq \mathcal{D}_{SRK}$,
  if and only if $a \geq \tfrac{1}{4}$ and  $b \in \mathbb{R} \setminus \{0\}$.
\end{Lem}
\textbf{Proof.} Inserting \eqref{Class-II-Gamma} and
\eqref{Class-II-Sigma-1-2} into \eqref{Stab-func} we have to prove
that
\begin{equation} \label{Proof-eqn-1}
  \begin{split}
    \hat{R}(\hat{h},k) = \frac{|a \hat{h}-1|^4 + \tfrac{1}{2} |k|^4 +
    |k|^2 \, |a \hat{h}-1|^2} {|a \hat{h}-1|^4 \, |\hat{h}-1|^2} < 1
    \, ,
  \end{split}
\end{equation}
because $\mathcal{D}_{SDE} \subseteq \mathcal{D}_{SRK}$ if
$\hat{R}(\hat{h},k) < 1$ for all $\hat{h}, k \in \mathbb{C}^2$ with
$2 \Re(\hat{h}) + |k|^2 <0$. Assuming that $\Re(\hat{h})<0$ and
$|k|^2 < -2 \Re(\hat{h})$, we prove that $\hat{R}(\hat{h},k)-1 < 0$
holds. Using this assumption, we get
\begin{equation} \label{Proof-eqn-2a}
    \hat{R}(\hat{h},k)-1 < \frac{\phi(\hat{h},a)}{|a \hat{h}-1|^4 \cdot |\hat{h}-1|^2}
\end{equation}
with
\begin{equation} \label{Proof-eqn-3}
    \phi(\hat{h},a) := |a \hat{h}-1|^4 + 2 \Re(\hat{h})^2 - 2 \Re(\hat{h}) \cdot |a \hat{h}-1|^2
    - |a \hat{h}-1|^4 \cdot |\hat{h}-1|^2 \, .
\end{equation}
Since the denominator in \eqref{Proof-eqn-2a} is positive, it is
sufficient to prove that $\phi(\hat{h},a) \leq 0$.
Considering \eqref{Proof-eqn-3} and collecting for the real part of $\hat{h}$
results in
\begin{equation} \label{Proof-eqn-4}
  \begin{split}
    \phi(\hat{h},a) = &-a^4 \Re(\hat{h})^6 + (4a^3+2a^4) \Re(\hat{h})^5 \\
    &+ (-8a^3-3a^4 \Im(\hat{h})^2
    -6a^2) \Re(\hat{h})^4 \\
    &+(4a^4 \Im(\hat{h})^2 + 10a^2+4a+8a^3 \Im(\hat{h})^2) \Re(\hat{h})^3 \\
    &+ (1-3a^4 \Im(\hat{h})^4 - 8a^3 \Im(\hat{h})^2 -8a^2 \Im(\hat{h})^2 -4a)
    \Re(\hat{h})^2 \\
    & + (4a \Im(\hat{h})^2 + 4a^3 \Im(\hat{h})^4 + 2a^4 \Im(\hat{h})^4
    + 2a^2 \Im(\hat{h})^2 ) \Re(\hat{h}) \\
    &- \Im(\hat{h})^2 -a^4 \Im(\hat{h})^6 -2a^2 \Im(\hat{h})^4 \, .
  \end{split}
\end{equation}
Due to our assumption $\Re(\hat{h})<0$, it is easy to see that
$\phi(\hat{h},a) \leq 0$ if $a \geq \tfrac{1}{4}$. Thus, we get
$A$-stability for $a \geq \tfrac{1}{4}$. \\ \\
As the final step, we prove that this bound for $a$ is also sharp.
Let us choose $\Im(\hat{h})=\Im(k)=0$.
For the proof, we consider
the boundary of the set $\mathcal{D}_{SDE}$ in the real case, which reduces to the half-line
$\partial \mathcal{D}_{SDE} := \{ (\hat{h}, \sqrt{-2\hat{h}}) : \hat{h} \in \, ]-\infty,0[ \, \}$.
Let $\psi(\hat{h},a) := \hat{R}(\hat{h}, \sqrt{-2 \hat{h}}) -1$. Then, we get from \eqref{Proof-eqn-1} that
\begin{equation} \label{Proof-eqn-7}
    \lim_{\hat{h} \to 0} \psi(\hat{h},a) = 0 \, .
\end{equation}
Now, the idea is to show that $\psi(\hat{h},a)$ is strictly decreasing
on the open set $\mathcal{S}_{\varepsilon} := \{ (\hat{h}, \sqrt{-2\hat{h}}) :
\hat{h} \in \, ]-\varepsilon,0[ \, \}$ for some
$\varepsilon=\varepsilon(a)>0$, i.e.\ $\frac{\partial \psi(\hat{h}, a)}{\partial \hat{h}}<0$.
Since $\hat{R}(\hat{h},k)$ is continuous on $\,]-\infty,0] \times \mathbb{R}$, for each
point $P \in \mathcal{S}_{\varepsilon}$ there exists an open ball $B_{\delta}(P)$ for some
$\delta>0$, such that $\hat{R}(\hat{h},k)>1$ on $B_{\delta}(P) \cap \mathcal{D}_{SDE}$,
i.e.\ the scheme is not $A$-stable. Thus, we consider
\begin{equation} \label{Proof-eqn-8}
  \begin{split}
  &\frac{\partial \psi(\hat{h}, a)}{\partial \hat{h}} = \frac{-2 \hat{h}}{(\hat{h} a -1)^5 (\hat{h}-1)^3} \\
  & \times (\hat{h}^4 a^5 -5 \hat{h}^3 a^4
  +(11 \hat{h}^2 -3 \hat{h}^3)a^3 + (7 \hat{h}^2 -11 \hat{h})a^2 +(4+4 \hat{h}^2 -7
  \hat{h})a-1) \, .
  \end{split}
\end{equation}
We distinguish the cases $a \in [0,\tfrac{1}{4}[\,$ and $a <0$. Let $-1 < \hat{h} < 0$.
First, we consider $a \in [0,\tfrac{1}{4}[ \,$. Then, using the estimates $\hat{h}^4<-\hat{h}$,
$-\hat{h}^3<-\hat{h}$ and $\hat{h}^2<-\hat{h}$ we obtain
\begin{equation} \label{Proof-eqn-9}
  \begin{split}
  \frac{\partial \psi(\hat{h}, a)}{\partial \hat{h}} < \frac{-2 \hat{h} (4a-1-\hat{h}(a^5+5a^4+14a^3+18a^2+11a))}
  {(\hat{h} a -1)^5 (\hat{h}-1)^3} \, .
  \end{split}
\end{equation}
Since $\frac{-2 \hat{h}}{(\hat{h} a -1)^5 (\hat{h}-1)^3}>0$ it obviously follows from
\eqref{Proof-eqn-9} that
$\frac{\partial \psi(\hat{h}, a)}{\partial \hat{h}}<0$ for $a=0$. Further, for $a \in \, ]0,\tfrac{1}{4}[\,$
we have $\frac{\partial \psi(\hat{h}, a)}{\partial \hat{h}}<0$ if
\begin{equation} \label{Proof-eqn-10}
  \max \left\{ \tfrac{4a-1}{a^5+5a^4+14a^3+18a^2+11a}, -1 \right\} < \hat{h} < 0 \, .
\end{equation}
For the case $a<0$, let $\max\{ \tfrac{1}{a}, -1\} < \hat{h} <0$.
Then, using the estimate $-\hat{h}^3 < \hat{h}^2$ and neglecting
some negative terms, we get from \eqref{Proof-eqn-8} that
\begin{equation} \label{Proof-eqn-11}
  \frac{\partial \psi(\hat{h}, a)}{\partial \hat{h}} < \frac{-2 \hat{h} (\hat{h}^2(5a^4+18a^2)-1)}
  {(\hat{h} a -1)^5 (\hat{h}-1)^3} \, .
\end{equation}
Since $\frac{-2 \hat{h}}{(\hat{h} a -1)^5 (\hat{h}-1)^3}>0$ it
follows that $\frac{\partial \psi(\hat{h}, a)}{\partial \hat{h}}<0$
if
\begin{equation} \label{Proof-eqn-12}
    \max \left\{ \tfrac{-1}{\sqrt{5a^4+18a^2}}, \tfrac{1}{a}, -1
    \right\} < \hat{h} < 0 \, ,
\end{equation}
which completes the proof. \hfill $\Box$ \\ \\
Considering class X in case of a diagonally drift implicit stiffly
accurate SRK method, we get with $a_1=a_2=a$ for some $a \in
\mathbb{R}$, $a_3=1$ and $a_4=0$ a family of SRK schemes
\eqref{SSRK} that need three stage-evaluations of the drift $f$ and
two stage-evaluations of the diffusion $g$ each step. Again, we need
$a \neq 0$ if the SRK method is applied to SDAEs.
\begin{Lem} \label{Lemma-A-stab-Class-X}
    The family of order 1.0 SRK schemes with coefficients
    \eqref{Sec:Opt-schemes-Class-X} in the case of $a_1=a_2$, $a_3=1$ and $a_4=0$ is
    $A$-stable for equation \eqref{Lin-stoch-test-eqn}, i.e.\
    $\mathcal{D}_{SDE} \subseteq \mathcal{D}_{SRK}$, if and only if $a_1
    \geq \tfrac{1}{4}$ and  $b \in \mathbb{R} \setminus \{0\}$.
\end{Lem}
\textbf{Proof.} The assertion follows from the fact that for
$a_1=a_2$, $a_3=1$ and $a_4=0$ the stability function \eqref{Stab-func} with
\eqref{Class-X-Gamma} and \eqref{Class-X-Sigma-1-2} coincides with
the stability function for the coefficients
\eqref{Sec:Opt-schemes-Class-II} of class II under the assumptions of
Lemma~\ref{Lemma-A-stab-Class-II}.
Therefore, the result follows from the proof of
Lemma~\ref{Lemma-A-stab-Class-II}. \hfill $\Box$
\begin{Bem}
    With the choice $a_1=a_2=a_3=1$ for both classes II and X, we get
    families of $A$-stable stiffly accurate singly diagonally
    drift-implicit SRK schemes.
    Therefore, the calculation of only one $LU$
    decomposition is needed each step if a simplified Newton method is applied
    to solve the implicit equations
    (see also \cite{Hai2}).
\end{Bem}
Next, we try to find within classes II and X some $A$-stable stiffly
accurate SRK schemes with a minimized number of stage-evaluations
for the drift function $f$ and the diffusion function $g$ needed
each step. Therefore, we analyse some stiffly accurate SRK methods
with an explicit first stage, i.e., we choose $a_1=0$ in the
following. These schemes can be applied to SDAEs as well, provided
that the sub-matrix $(A_{ij})_{2 \leq i,j \leq s}$ is nonsingular
\cite{KKR12}.
\begin{figure}
\begin{center}
\includegraphics[width=6.8cm]{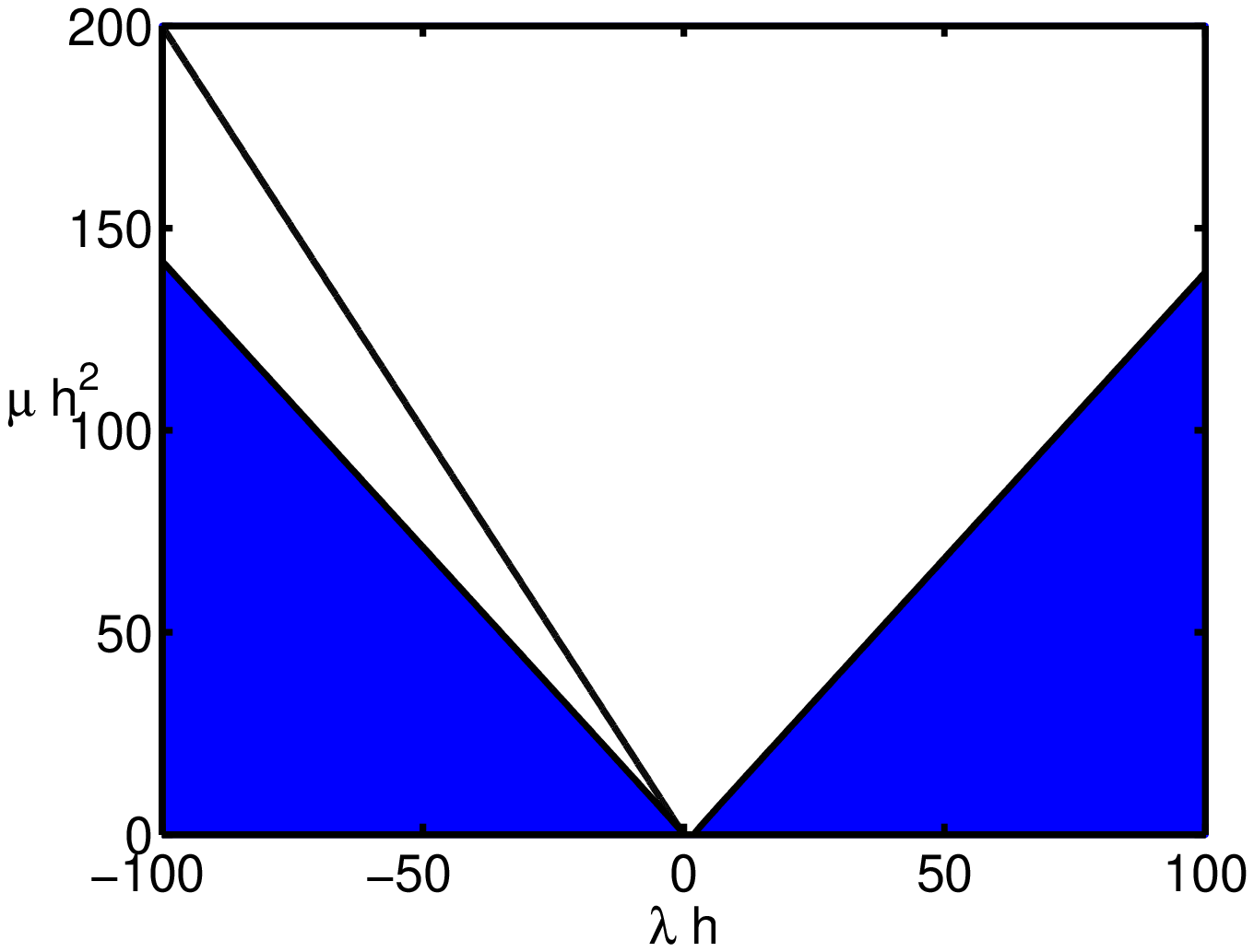}
\includegraphics[width=6.8cm]{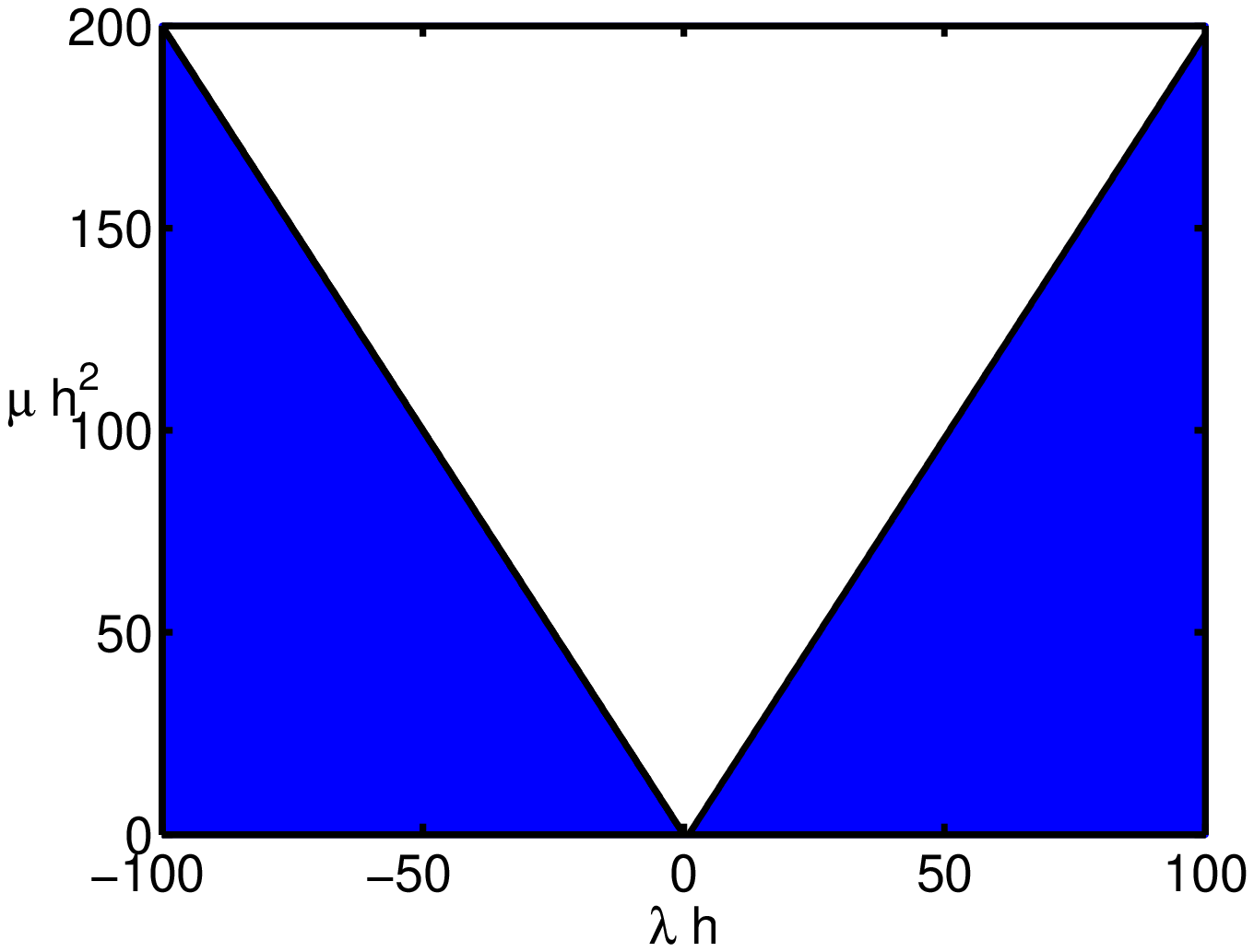}
\includegraphics[width=6.8cm]{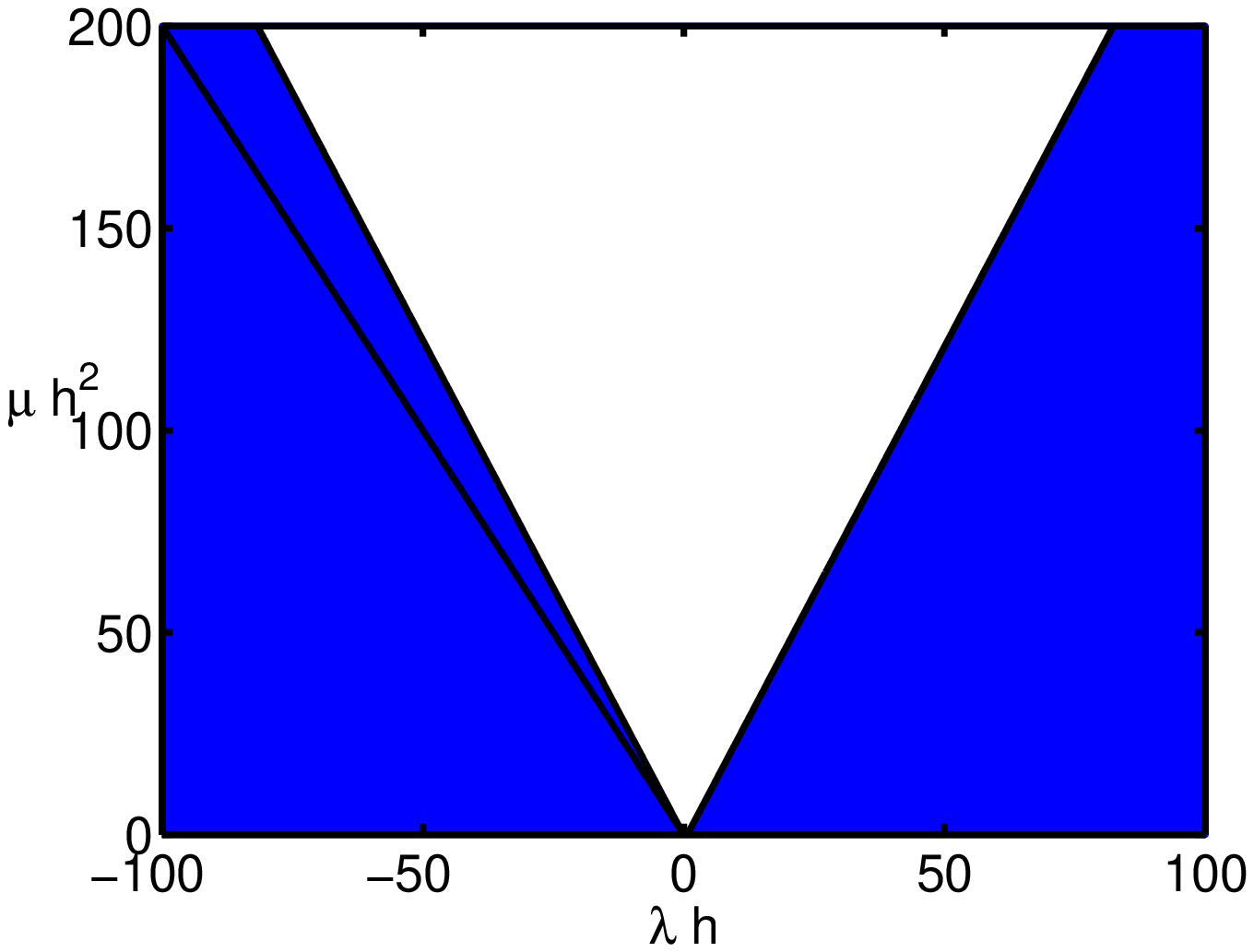}
\includegraphics[width=6.8cm]{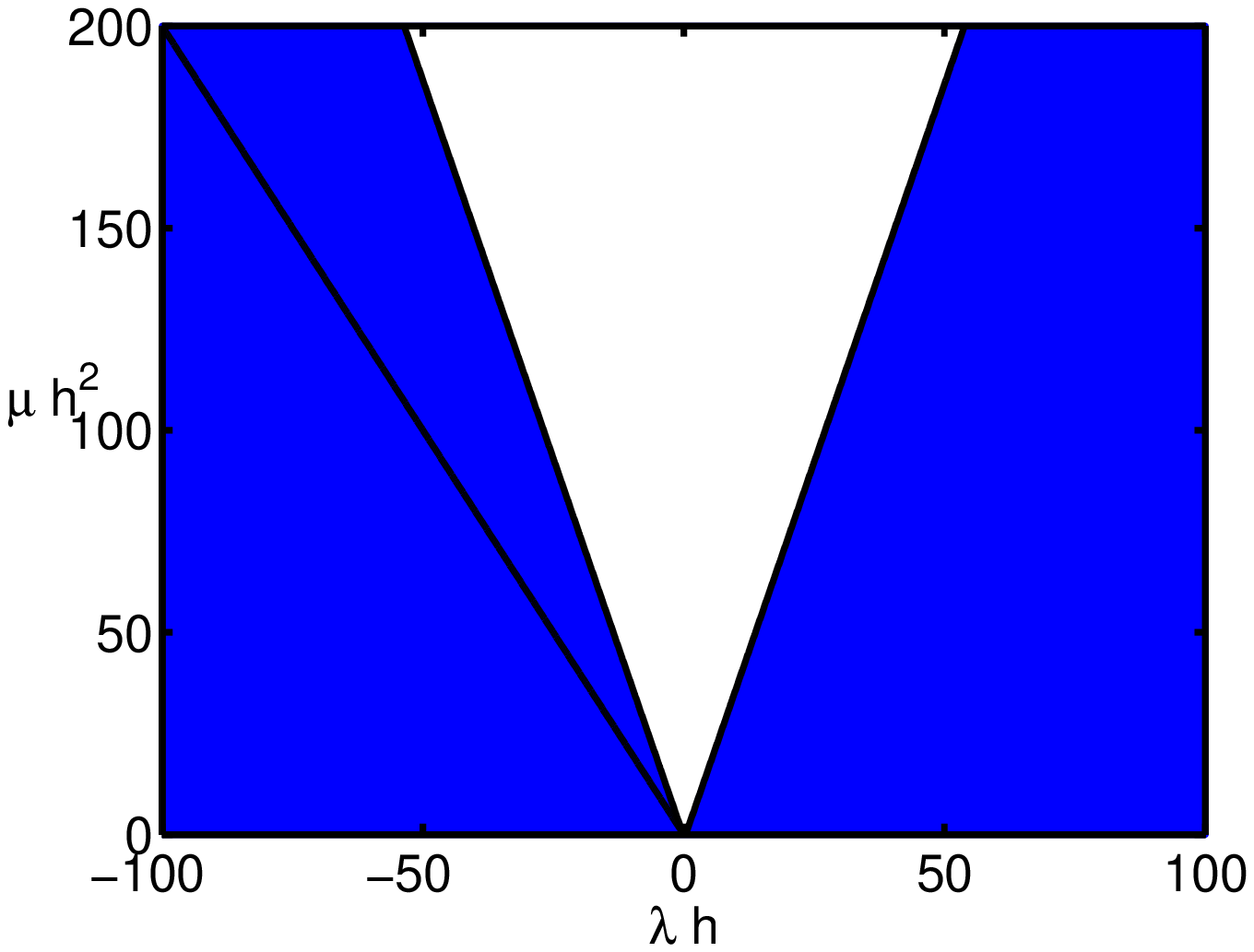}
\caption{Mean-square stability region for class II and class X
applicable only to SDEs and not to SDAEs with $a_1=a_2=a_4=0$, $b=1$
and with $a_3=1$ (not $A$-stable), $a_3=\frac{3}{2}$, in the lower
figures with $a_3=2$ and $a_3=4$, respectively.}
\label{MS-Bild-class-II-X-2-SDE}
\end{center}
\end{figure}
\begin{figure}
\begin{center}
\includegraphics[width=6.8cm]{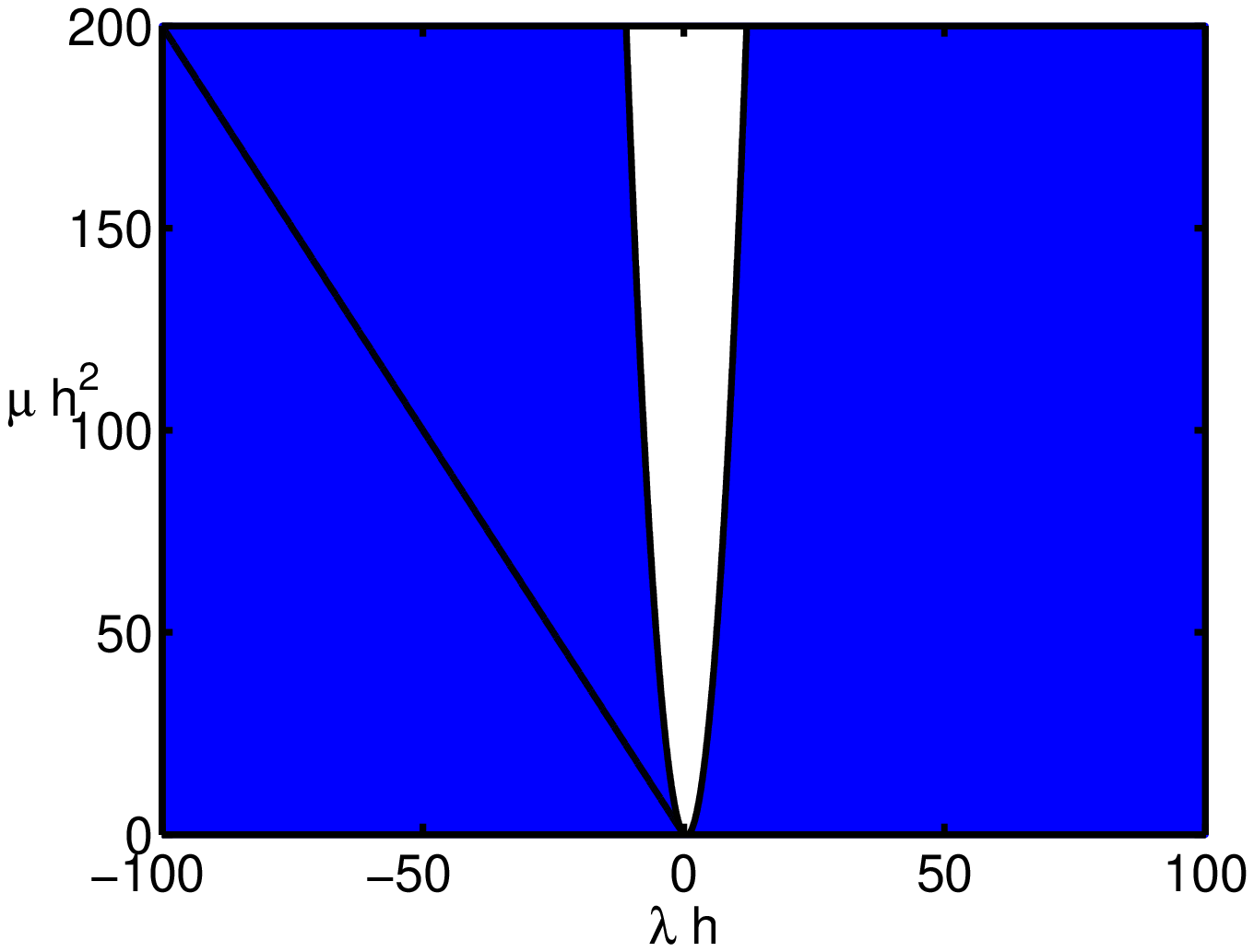}
\includegraphics[width=6.8cm]{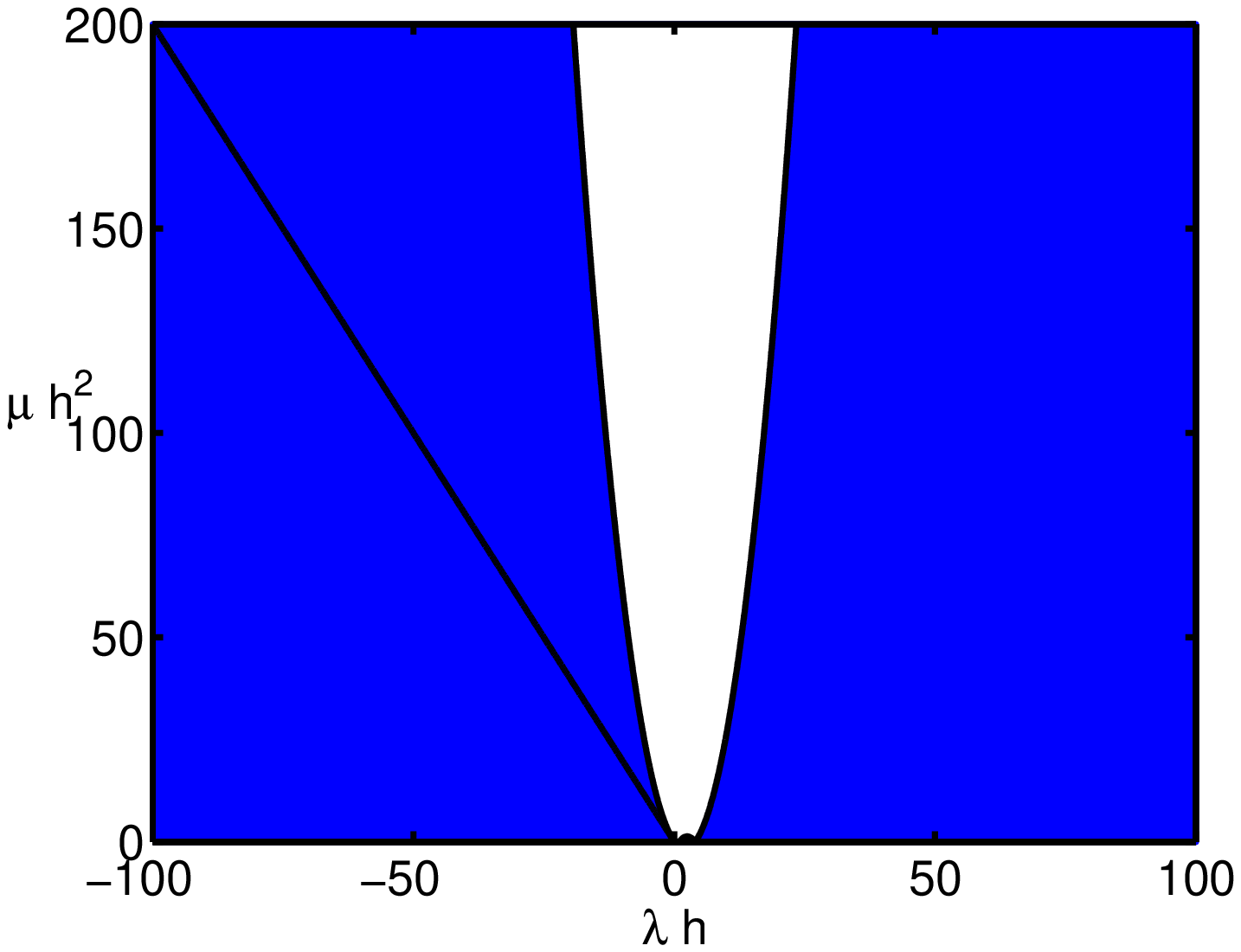}
\includegraphics[width=6.8cm]{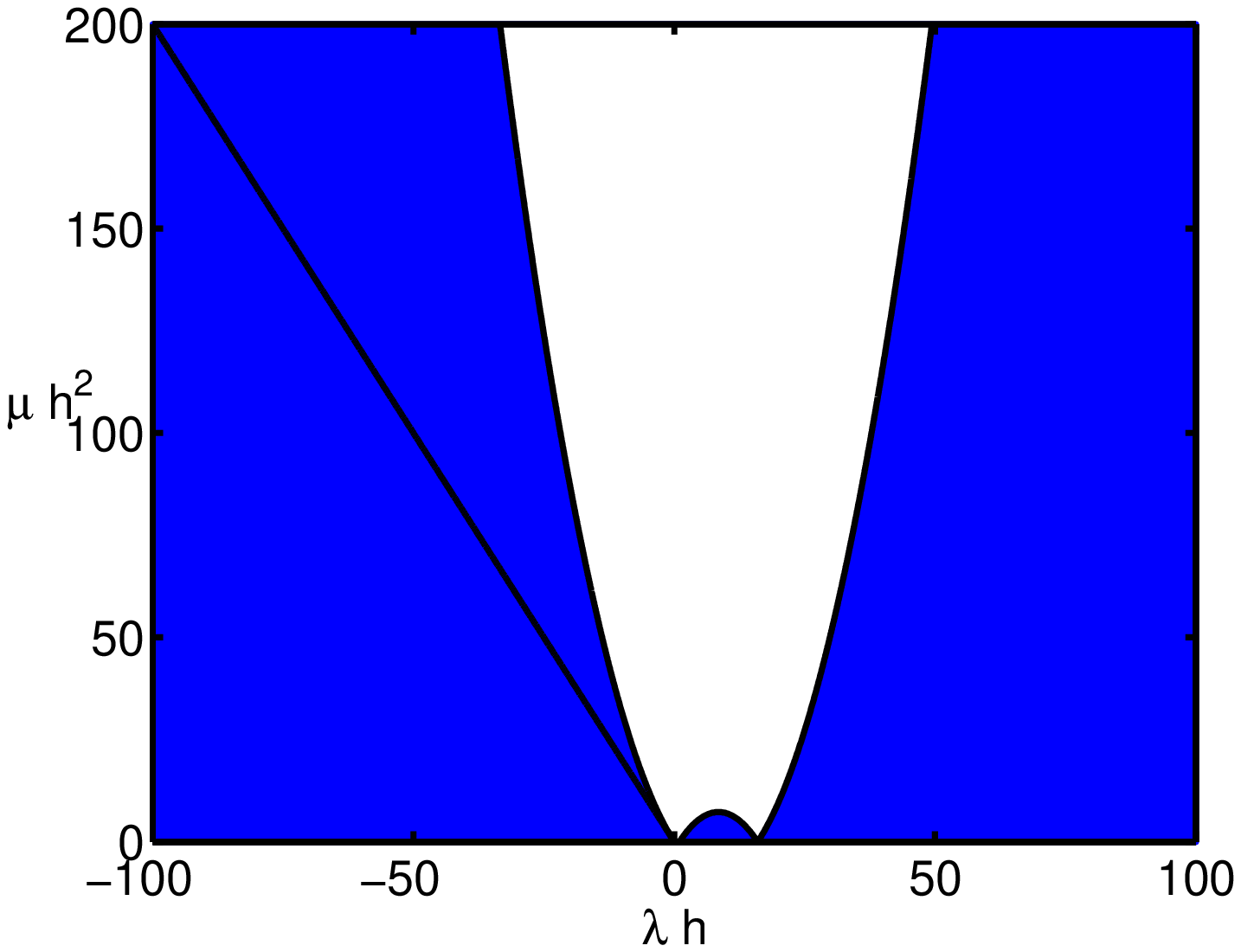}
\includegraphics[width=6.8cm]{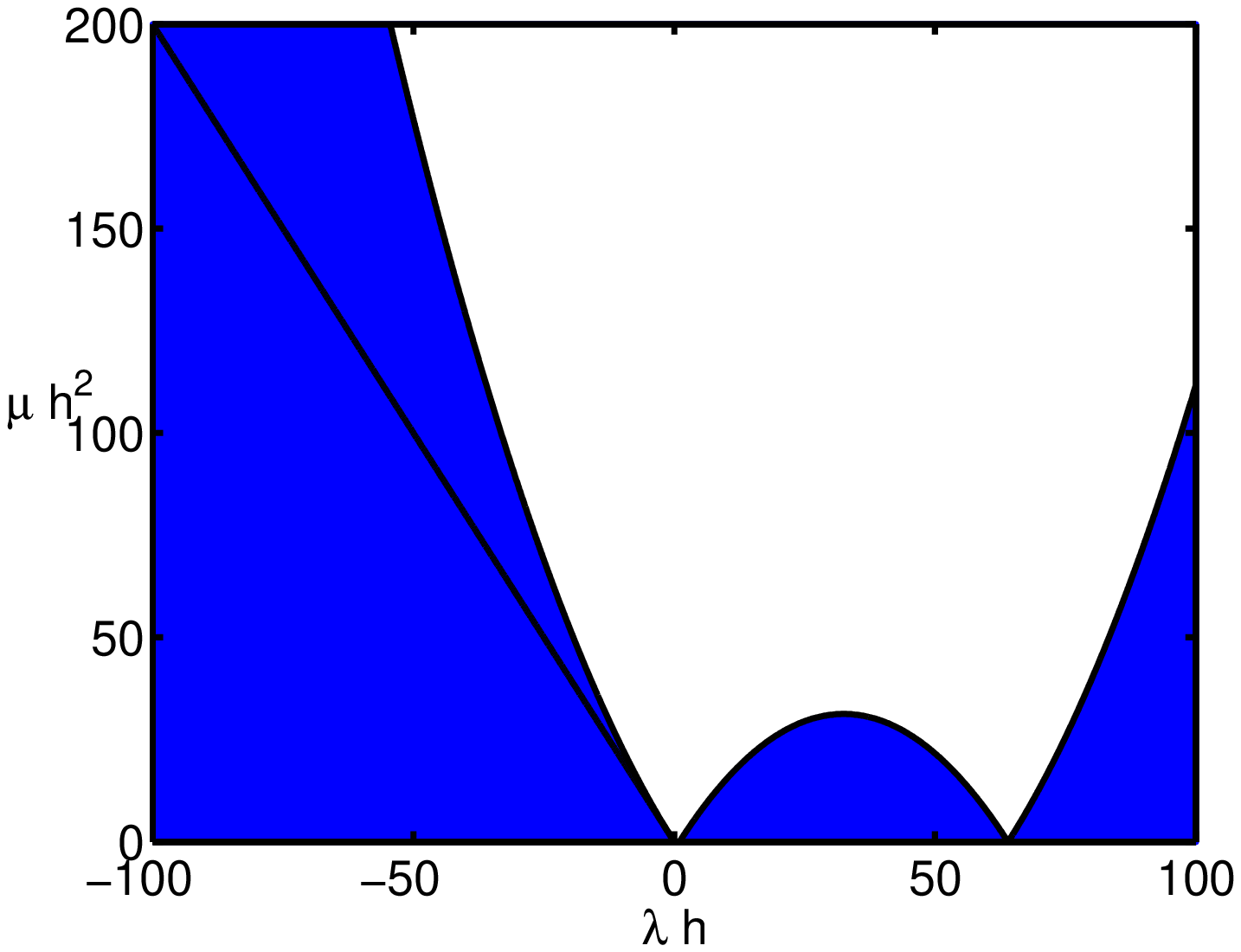}
\caption{Mean-square stability region for class II and class X
applicable to SDAEs with $a_1=0$, $a_3=\frac{3}{2}$, $a_4=-a_2$,
$b=1$ and with $a_2=\frac{3}{2}$, $a_2=\frac{1}{4}$, in the lower
figures with $a_2=\frac{1}{16}$ and $a_2=\frac{1}{64}$,
respectively.} \label{MS-Bild-class-II-X-2-SDAE}
\end{center}
\end{figure}
\begin{Lem} \label{Lemma-A-stab-Class-II-opt}
  The family of order 1.0 SRK schemes with coefficients
  \eqref{Sec:Opt-schemes-Class-II}
  is $A$-stable for equation \eqref{Lin-stoch-test-eqn} in case of $a_1=0$,
  i.e.\ it holds $\mathcal{D}_{SDE} \subseteq \mathcal{D}_{SRK}$,
  if $a_2 \geq 0$, $a_3 \geq \tfrac{3}{2}$ and  $b \in \mathbb{R} \setminus \{0\}$.
\end{Lem}
\textbf{Proof.} First, calculate the stability function
$\hat{R}(\hat{h},k)$ from \eqref{Stab-func}. As a result of this, we
have to prove that
\begin{equation} \label{Proof-3-eqn-1}
  \begin{split}
    \hat{R}(\hat{h},k) = \frac{\tfrac{1}{2} |k|^4 + |a_2
    \hat{h}-1|^2 \,
    (|k|^2 + |\hat{h}|^2 (1-a_3)^2 +
    (\hat{h}+\overline{\hat{h}})(1-a_3) + 1)}
    {|a_2 \hat{h}-1|^2 \, |a_3 \hat{h}-1|^2} < 1 \, ,
  \end{split}
\end{equation}
for all $\hat{h}, k \in \mathbb{C}^2$ with $2 \Re(\hat{h}) + |k|^2
<0$. Therefore, we assume that $\Re(\hat{h})<0$ and $|k|^2 < -2
\Re(\hat{h})$ and we prove that $\hat{R}(\hat{h},k)-1 < 0$ is
fulfilled under the assumptions of
Lemma~\ref{Lemma-A-stab-Class-II-opt}. Using these assumptions, we
get
\begin{equation} \label{Proof-3-eqn-2a}
    \begin{split}
    \hat{R}(\hat{h},k)-1 &= \frac{\frac{1}{2} |k|^4 + |a_2 \hat{h}-1|^2
    (|k|^2 + |\hat{h}|^2 -2a_3 |\hat{h}|^2 + 2 \Re(\hat{h}))}{|a_2
    \hat{h}-1|^2 \cdot |a_3 \hat{h}-1|^2} \\
    &< \frac{\phi(\hat{h},a_2,a_3)}{|a_2
    \hat{h}-1|^2 \cdot |a_3 \hat{h}-1|^2}
    \end{split}
\end{equation}
where
\begin{equation} \label{Proof-3-eqn-3}
    \phi(\hat{h},a_2,a_3) := 2 \Re(\hat{h})^2 + (\Re(\hat{h})^2 + \Im(\hat{h})^2) (1-2a_3)
    |a_2 \hat{h}-1|^2 .
\end{equation}
The denominator in \eqref{Proof-3-eqn-2a} is positive, thus it is
sufficient to prove $\phi(\hat{h},a) \leq 0$.
Collecting for the real part of $\hat{h}$ in \eqref{Proof-3-eqn-3},
we get
\begin{equation} \label{Proof-3-eqn-4}
  \begin{split}
    \phi(\hat{h},a_2,a_3) = &\, a_2^2(1-2a_3) \Re(\hat{h})^4
    + (-2a_2(1-2a_3)) \Re(\hat{h})^3 \\
    &+ (3 -2a_3 +2(1-2a_3) a_2^2 \Im(\hat{h})^2) \Re(\hat{h})^2 \\
    &-2 a_2(1-2a_3) \Im(\hat{h})^2 \Re(\hat{h}) \\
    &+(1-2a_3) \Im(\hat{h})^2 + (1-2a_3) a_2^2 \Im(\hat{h})^4 \, .
  \end{split}
\end{equation}
Due to our assumption $\Re(\hat{h})<0$, it is easy to see that
$\phi(\hat{h},a_2,a_3) \leq 0$ if $3-2a_3 \leq 0$, $1-2a_3 \leq 0$
and $a_2 \geq 0$. Thus, we need $a_2 \geq 0$ and $a_3 \geq
\tfrac{3}{2}$ for $A$-stability. \hfill $\Box$
\\ \\
%
%
%
%
Considering now class X with an explicit first stage, that is in
case of $a_1=0$, then similar results can be obtained if the
simplifying assumption $a_4=-a_2$ is fulfilled.
\begin{Lem} \label{Lemma-A-stab-Class-X-opt}
  The family of order 1.0 SRK schemes with coefficients
  \eqref{Sec:Opt-schemes-Class-X}
  is $A$-stable for equation \eqref{Lin-stoch-test-eqn} in case of $a_1=0$ and $a_4=-a_2$,
  i.e.\ it holds $\mathcal{D}_{SDE} \subseteq \mathcal{D}_{SRK}$,
  if $a_2 \geq 0$, $a_3 \geq \tfrac{3}{2}$ and  $b \in \mathbb{R} \setminus \{0\}$.
\end{Lem}
\textbf{Proof.} We calculate the stability function
$\hat{R}(\hat{h},k)$ from \eqref{Stab-func} in the case of $a_1=0$. Therefore, we
have to prove that
\begin{equation} \label{Proof-4-eqn-1}
  \begin{split}
    \hat{R}(\hat{h},k) &= \frac{\tfrac{1}{2} |k|^4 + |1-a_2 \hat{h}|^2 \,
    (|k|^2 + |\hat{h}|^2 (1-2a_3) +
    \hat{h}+\overline{\hat{h}} + |1-a_3 \hat{h}|^2)}
    {|a_2 \hat{h}-1|^2 \, |a_3 \hat{h}-1|^2} \\
    &+ \frac{\frac{1}{2b} |k|^2 (a_2+a_4) (\overline{k} \hat{h} + k \overline{\hat{h}})
    + \frac{1}{2 b^2} |k|^2 |\hat{h}|^2 (a_2+a_4)^2}{|a_2 \hat{h} -1|^2 \, |a_3 \hat{h} -1|^2} < 1 \, ,
  \end{split}
\end{equation}
for all $\hat{h}, k \in \mathbb{C}^2$ with $2 \Re(\hat{h}) + |k|^2
<0$. Thus, we assume that $\Re(\hat{h})<0$ and $|k|^2 < -2
\Re(\hat{h})$ and we prove that $\hat{R}(\hat{h},k)-1 < 0$ is
fulfilled under the assumptions of
Lemma~\ref{Lemma-A-stab-Class-X-opt}. Using the assumption that additionally $a_4=-a_2$,
we calculate that
\begin{equation} \label{Proof-4-eqn-2a}
    \begin{split}
    \hat{R}(\hat{h},k)-1 &= \frac{\frac{1}{2} |k|^4 + |a_2 \hat{h}-1|^2
    (|k|^2 + |\hat{h}|^2 -2a_3 |\hat{h}|^2 + 2 \Re(\hat{h}))}{|a_2
    \hat{h}-1|^2 \cdot |a_3 \hat{h}-1|^2} . 
    \end{split}
\end{equation}
Here, it turns out that \eqref{Proof-4-eqn-2a} is the same as \eqref{Proof-3-eqn-2a}, i.e.,
the rest of the proof is the same as in the proof of Lemma~\ref{Lemma-A-stab-Class-II-opt}.
\hfill $\Box$
\\ \\
Figure~\ref{MS-Bild-class-II-X-2-SDE} shows the region of
mean-square stability for classes II and X with $a_1=a_2=a_4=0$,
$b=1$ and for $a_3 \in \{1, \tfrac{3}{2}, 2, 4\}$. We point out that
the SRK schemes with $a_1=a_2=0$ can be applied to SDEs but not to
SDAEs. Here, the SRK schemes SADIRK12II and SADIRK12X for SDEs with
coefficients $a_1=a_2=a_4=0$, $b=1$ and $a_3=\tfrac{3}{2}$ for class
II and X, respectively, cover the region of stability for the SDE
best. However, for SDAEs we need $a_2 \neq 0$ and $a_3 \neq 0$ if
$a_1=0$ holds. Some regions of mean-square stability for the
coefficients $a_1=0$, $a_3=\tfrac{3}{2}$, $b=1$ and $a_2 \in
\{\tfrac{3}{2}, \tfrac{1}{4}, \tfrac{1}{16}, \tfrac{1}{64}\}$ are
given in Figure~\ref{MS-Bild-class-II-X-2-SDAE}. Here, it can be
observed that the region of mean-square stability fits better to the
region of the test equation the smaller the values of $a_2$ are.
\begin{Bem}
    If we choose $a_1=a_2=0$ and $a_3 \geq \tfrac{3}{2}$ for class II
    or for class X with additionally $a_4=0$,
    then only one stage-evaluation of the drift $f$ and two stage-evaluations of the
    diffusion $g$ are necessary each step for the $A$-stable stiffly
    accurate SRK scheme \ref{SSRK} applicable to SDEs. This is due to the FSAL (first
    same as last) property and due to an explicit first stage, see
    also e.~g.\ \cite{Hai2}. Further, we get a family of
    drift-implicit SRK schemes that need only
    one implicit equation to be solved each step. However, to apply the
    stiffly accurate SRK schemes to SDAEs, we need $a_2 \neq 0$ and
    thus two stage-evaluations of the drift $f$ and two stage-evaluations of the
    diffusion $g$ due to FSAL. Especially, in case of $a_2=a_3$ only
    one $LU$ decomposition has to be calculated each step if a
    simplified Newton method is applied to solve the implicit equations
    (see \cite{Hai2}).
\end{Bem}
\section{Conclusions}
We have calculated a classification of the set of solutions for the
order conditions of stiffly accurate strong order 0.5 and order 1.0
SRK methods for SDAEs with a scalar driving Wiener process
introduced in \cite{KKR12}. As the main advantages of the considered
SRK method compared to well known schemes, no projectors and no
pseudo-inverses have to be calculated and the considered SRK methods
are derivative-free what makes them easy to be implemented. Based on
this classification, a mean-square stability analysis is carried out
for the two classes II and X. These two classes allow to minimize
the computational costs in the sense that a minimum number of
stage-evaluations are needed as well as a minimum number of implicit
equations that have to be solved each step. Further, the two classes
represent both cases where $B^{(2)} = 0$ and $B^{(2)} \neq 0$, i.e.,
where the random variables $I_{(1,1),n}$ do not appear and do appear
explicitly within the scheme, respectively. For both classes II and
X, conditions for the coefficients such that the SRK method is
$A$-stable in the mean-square sense are proved for diagonally
drift-implicit schemes as well as for schemes with an explicit first
stage. Especially, a family of $A$-stable stiffly accurate
drift-implicit order 1.0 SRK schemes for SDEs has been found that
needs only one stage-evaluation of the drift function $f$, two
stage-evaluations of the diffusion function $g$ and one implicit
equation to be solved each step. However, for the SDAE case at least
two stage-evaluations of the drift $f$ and two stage-evaluations of
the diffusion $g$ are needed for an $A$-stable SRK method. For
future research it would be interesting to analyse not only
mean-square stability, but maybe to find some further concepts of
stability that are of importance especially for SDAEs.
\bibliographystyle{abbrv}
\bibliography{SDAE}
\end{document}